\pdfoutput=1

\RequirePackage{fix-cm}

\documentclass[leqno, fleqn, centertags, 12pt]{article}

\usepackage{latexsym}
\usepackage{amsmath}
\usepackage{amssymb}
\usepackage{verbatim}

\usepackage[T1]{fontenc}

\usepackage[upright, widespace]{fourier}

\usepackage{fullpage}

\makeatletter
\renewcommand{\@makefntext}[1]{\vspace*{0.5ex}\parindent=0em
\hspace*{-0.4em}
\hbox to 0.4em{\hss\@makefnmark}\hspace*{0.4em}{#1}
}
\makeatother

\newcounter{mysectionnumber}
\newcommand{\mysection}[2]{
\refstepcounter{mysectionnumber}
\section*{ \textnormal{{\themysectionnumber.}\oss {#1}}}\label{#2}}

\newcommand{\mynonumbersection}[1]{
\vspace{-0.0ex}
\section*{{}\hspace*{0.00em}$\phantom{1.}$\textnormal{{#1}}}}

\newcommand{\myit}[1]{\textbf{\textit{#1}}\hspace{0.0em}}

\newcounter{myparnum}
\renewcommand{\themyparnum}{\arabic{myparnum}}
\newcommand{\mypar}[2]{\refstepcounter{myparnum}{\vspace{\medskipamount}\textbf{{\themyparnum. #1}\label{#2}}\hspace{0.5em}}}

\newcounter{mylemmanum}[myparnum]
\renewcommand{\themylemmanum}{\themyparnum.\arabic{mylemmanum}}
\newcommand{\mylemma}[2]{\refstepcounter{mylemmanum}{\vspace{\medskipamount}\textbf{\textit{\themylemmanum. #1}\label{#2}}\hspace{0.5em}}}

\newcommand{\myuppar}[1]{\vspace{\medskipamount}\textbf{#1}\hspace*{0.5em}}

\newcommand{\proof}{\vspace{\medskipamount}{\textbf{{\emph{Proof}.}}\hspace*{1em}}}
\newcommand{\prooftitle}[1]{\vspace{\medskipamount}{\textbf{{\emph{#1}.}}\hspace*{1em}}}
\newcommand{\subproof}{\vspace{\medskipamount}{{\emph{Proof}.}\hspace*{1em}}}
\newcommand{\eproof}{ $\blacksquare$}
\newcommand{\esubproof}{ $\square$}

\newcommand{\dis}{\displaystyle}

\def\sss{\hspace{0.05em}\ }
\def\dss{\hspace{0.1em}\ }
\def\trs{\hspace{0.15em}\ }
\def\qss{\hspace{0.2em}\ }
\def\pss{\hspace{0.3em}\ }
\def\oss{\hspace{0.4em}\ }

\def\halfff{\hspace*{0.025em}}
\def\fff{\hspace*{0.05em}}
\def\dff{\hspace*{0.1em}}
\def\trf{\hspace*{0.15em}}
\def\qff{\hspace*{0.2em}}
\def\pff{\hspace*{0.3em}}
\def\off{\hspace*{0.4em}}

\newcommand{\hnsp}{\hspace*{-0.05em}}
\newcommand{\nsp}{\hspace*{-0.1em}}
\newcommand{\nnsp}{\hspace*{-0.15em}}

\newcommand{\dnsp}{\hspace*{-0.2em}}

\renewcommand{\leq}{\leqslant}
\renewcommand{\geq}{\geqslant}

\newcommand{\rrr}{\mathbf{R}}

\newcommand{\num}[1]{|\qff #1 \qff|}

\newcommand{\toto}{\longrightarrow}
\newcommand{\ttoo}{\hspace*{0.2em}\longrightarrow\hspace*{0.2em}}

\begin{document}

\setlength{\baselineskip}{12pt plus 0pt minus 0pt}
\setlength{\parskip}{12pt plus 0pt minus 0pt}
\setlength{\abovedisplayskip}{12pt plus 0pt minus 0pt}
\setlength{\belowdisplayskip}{12pt plus 0pt minus 0pt}

\newskip\smallskipamount \smallskipamount=3pt plus 0pt minus 0pt
\newskip\medskipamount   \medskipamount  =6pt plus 0pt minus 0pt
\newskip\bigskipamount   \bigskipamount =12pt plus 0pt minus 0pt

\author{Nikolai\qss V.\qss Ivanov}
\title{Beyond\pss Sperner's\pss lemma}
\date{}

\footnotetext{\hspace*{-0.65em}\copyright\oss 
Nikolai\qss V.\qss Ivanov,\oss 2017,\qss 2022\qss
(revised\sss version).\oss
The author\dss is\dss grateful\dss to\qss V.{\dff}I.\dss Danilov\dss
for attracting\dss his attention\dss to\sss the works of\qss Scarf\qss
and\dss very\sss stimulating\sss correspondence.}

\maketitle

\vspace*{6ex}

{\renewcommand{\baselinestretch}{1}
\selectfont

\myit{\hspace*{0em}\large Contents}\vspace*{1.5ex} \\ 
\hbox to 0.8\textwidth{\myit{Preface} \hfil 1}\hspace*{0.5em} \vspace*{1.5ex}\\
\hbox to 0.8\textwidth{\myit{\phantom{1}1.}\hspace*{0.5em} Scarf's\qss combinatorial\dss
theorem\dss in\sss an abstract\sss setting \hfil 3}\hspace*{0.5em} \vspace*{0.25ex}\\
\hbox to 0.8\textwidth{\myit{\phantom{1}2.}\hspace*{0.5em} Brouwer's\qss fixed\dss point\dss theorem \hfil 9}\hspace*{0.5em} \vspace*{0.25ex}\\
\hbox to 0.8\textwidth{\myit{\phantom{1}3.}\hspace*{0.5em} Scarf's\qss proof\dss
and\dss its\sss versions  \hfil 11}\hspace*{0.5em}\vspace*{1.5ex}\\
\hbox to 0.8\textwidth{\myit{References}\hspace*{0.5em}\hfil 14}\hspace*{0.5em}    

}

\renewcommand{\baselinestretch}{1}
\selectfont

\mynonumbersection{Preface}

\vspace*{6pt} 
In\qss 1967\qss H.\dss Scarf\pss \cite{sc2}\qss suggested a new\dss
proof\dss of\qss Brouwer's\qss fixed\dss point\dss theorem.\oss
In outline,\oss his proof\trs is\dss similar\sss to\sss the 
well\dss known\sss proof\dss based on\qss Sperner's\trs lemma and\dss
Knaster--Kuratowski--Mazurkiewicz\qss \cite{kkm}\qss argument\halfff.\oss
But\dss instead of\trs triangulations\dss Scarf\trs uses\qss (sufficiently\dss dense)\qss 
finite subsets of\dss a simplex\halfff.\oss
This\sss forces a modification of\dss the\dss
KKM\qss argument\halfff:\oss
the inequalities in\dss the coloring\dss rule must\dss be reversed.\oss
Also,\oss instead of\dss using\sss the fundamental\dss fact\sss that\sss an odd\sss number cannot\dss
be equal\dss zero,\oss Scarf\trs use a path-following\sss algorithm\dss motivated\sss by\dss
the\sss game\sss theory and\dss the\sss linear programming.\oss
From\dss Scarf's\dss point\sss of\dss view,\oss the algorithmic nature of\dss
his proof\dss was especially\sss important.\oss
But,\oss as explained\sss in\qss \cite{i1},\oss the proof\dss based on\dss Sperner's\trs lemma
can be easily\dss turned\dss to a path-following algorithm.\oss
Moreover,\oss path-following\sss arguments we used already\dss in\qss 1929\qss by\trs Hurewicz\qss \cite{h}.\oss

The most\sss striking\sss feature of\qss
Scarf's\trs proof\qss is\dss a\sss combinatorial\dss theorem,\oss 
namely,\oss an analogue of\qss Sperner's\trs lemma\sss dealing\sss with\sss
finite subsets of\dss a simplex\sss instead of\trs triangulations.\oss 
See\qss \cite{sc2},\oss Theorem\qss 1.\oss
The origins of\qss this\sss theorem and\dss its\sss proof\trs also belong\dss
to\sss the\sss linear\dss programming\dss and\dss the game\sss theory\halfff.\oss
Actually\halfff,\oss Scarf's\qss proof\dss of\qss Brouwer's\qss theorem\dss is\dss
a\sss byproduct\sss of\pss Scarf's\qss fundamental\dss paper\qss \cite{sc1}\qss
in\dss game\sss theory\halfff.\oss 
A\dss topological\dss interpretation of\qss Scarf's\qss combinatorial\dss theorem 
was suggested\dss by\qss H.\dss Kuhn.\oss 
See\qss \cite{ku},\oss {\sc Scarf's\dss Theorem}.\vspace{-0.55pt} 

Scarf's\qss combinatorial\dss theorem\dss is\dss a\sss truly\sss combinatorial\dss
result\halfff,\oss
but\dss the\sss style of\dss his papers\dss is\dss very\dss geometric,\pss
and\dss he always puts forward\sss the geometric underpinnings of\dss combinatorial\sss
arguments.\oss
For example,\oss Scarf\qss orders subsets of\dss a simplex\dss by\dss the values of\dss
one of\qss barycentric coordinates,\oss but\dss
instead of\trs taking\dss the minimal\sss element\dss of\dss a\sss set\dss
he\sss moves a face of\dss a simplex\dss parallely\dss to\sss itself\qss
(see\qss \cite{sc2},\oss the proof\dss of\qss Lemma\qss 1,\pss for example).\oss\vspace{-0.55pt}

Recently\qss H.\qss Petri\dss and\dss M.\qss Voorneveld\qss \cite{pv}\qss
published\sss an almost\trs geometry-free\sss version of\qss 
Scarf's\qss proof\dss of\qss Brouwer's\qss theorem.\oss
Eliminating\dss geometry\dss is\dss an explicitly\sss stated\dss goal\sss
of\pss \cite{pv}.\oss
See\qss \cite{pv},\oss Concluding\dss remarks.\oss
This\dss is\dss achieved\dss not\dss without\sss a cost\halfff,\oss
Petri--Voorneveld\qss version\sss of\pss Scarf's\qss proof\qss is\trs fairly\dss tricky\halfff.\oss
A minimal\sss amount\sss of\trs geometry\dss is\dss inevitable
since\qss Brou\-wer's\qss theorem\dss is\dss about\sss 
a simplex\halfff.\oss 
But\dss the geometry\sss of\dss simplices\dss is\dss
eliminated\dss by\dss dealing\dss with\sss
a\sss family\sss of\dss arbitrary\dss linear\sss orders 
on a\sss finite set\dss instead of\trs
the linear orders induced\dss by\dss barycentric coordinates
on a subset\sss of\dss a simplex\halfff.\oss
Scarf\qss himself\trs had applied\dss his\sss methods\sss to other orders
and\sss was aware\sss that\dss they\dss work for arbitrary\sss orders.\oss
See\qss \cite{sc3},\oss Chapter\qss 6\qss and,\pss especially\halfff,\pss
the discussion at\dss the bottom of\pss p.\trs 146.\oss
Rewriting\qss Scarf's\qss arguments\dss in\sss such an abstract\dss setting\dss
is\dss only\sss natural\sss and\sss makes\sss their 
combinatorial\sss nature\sss transparent\halfff.\oss\vspace{-0.55pt}

But\qss Petri\dss and\dss Voorneveld\trs went\dss further and eliminated\sss
also only\dss implicitly\dss geometric aspects of\trs the proof\halfff.\oss
Like\qss Sperner\qss \cite{s},\oss Petri\dss and\dss Voorneveld\qss \cite{pv}\qss start\dss with
a coloring\sss of\dss a finite set\halfff.\oss
In\qss Sperner's\qss context\dss this finite set\dss is\dss the set\sss
of\trs vertices of\dss a\sss triangulation,\oss
and\dss this\sss triangulation\dss is\dss present\dss 
before any\sss coloring\trs is\dss introduced\sss and\sss
even\dss before\sss
the proof\trs begins.\oss
A similar\dss structure on an arbitrary\sss finite set\sss of\dss points\sss in a\sss simplex\qss
(or\sss a\sss finite set\dss with a\sss family\sss of\trs
linear orders)\qss
was discovered\dss by\qss Scarf\pss \cite{sc1},\oss \cite{sc2}.\oss
It\trs is\trs the collection of\dss so-called\qss \emph{primitive\sss sets},\oss 
and\dss it\dss was\sss interpreted\dss by\qss Kuhn\qss \cite{ku}\qss
as\sss the structure of\dss a simplicial\sss complex\halfff.\oss
See\qss Section\qss \ref{comparing}\qss below.\oss
Petri\dss and\dss Voorneveld\trs use only\qss
(a version of\dff)\qss primitive sets
closely\dss related\dss to\sss the coloring\dss in\sss question,\pss
and\dss this renders\dss the\dss Scarf--Kuhn\dss structure\sss invisible.\oss\vspace{-0.55pt}

The present\dss paper\dss originated\dss in\dss a\sss text\dss
written\sss simply\dss to verify\qss Petri--Voorneveld\qss proof\halfff.\oss
After\dss the verification,\oss the proof\dss was rearranged\dss in order\dss
to stress similarities with\qss Sper\-ner's\qss proof\dss of\trs his\sss lemma.\oss
The rearranged\dss proof\trs turned out\dss to be equally\dss well\sss
suited\dss for comparing\dss with\qss Scarf's\oss one.\oss
After\sss such a comparison\dss it\dss became clear\dss
that\dss the\dss Scarf--Kuhn\dss structure should\dss be separated\dss from colorings.\oss\vspace{-0.55pt}

The main\dss goal\dss of\trs this\dss paper\dss is\dss to present\dss
the resulting\dss proof\dss
of\pss Scarf's\qss combinatorial\dss theorem\sss in\dss the
setting of\dss a\sss family\sss of\dss orders on a\sss finite set\halfff.\oss
This\dss is\dss done in\dss Section\qss \ref{scarf}.\oss
Section\qss \ref{brouwer}\qss is\dss devoted\dss to a deduction of\qss
Brouwer's\qss fixed\dss point\dss theorem\dss from\qss Scarf's\qss theorem.\oss 
This deduction\dss follows\qss
Petri\dss and\dss Voorneveld\qss \cite{pv}\qss version.\oss
Since\sss the geometry\dss is\dss avoided,\oss
their\sss version\dss is\dss a\dss little more complicated\dss than\qss Scarf's\qss one.\oss
Finally\halfff,\oss Section\qss \ref{comparing}\qss explains how\dss the proofs\dss
in\qss \cite{pv}\qss and\dss in\dss
the present\dss paper are related\dss to\qss Scarf's\qss proofs.\oss
In a\sss related\dss paper\qss \cite{i2}\qss Scarf's\qss results\trs are\dss
discussed\dss from\dss the point\sss of\dss view of\trs the combinatorial\dss
topology\halfff,\oss 
and,\oss in\dss particular\halfff,\oss some ideas\sss of\pss Kuhn\qss \cite{ku}\qss 
are refined.\oss
See\qss \cite{i2},\oss Sections\qss 1 -- 3.\oss

\newpage
\mysection{Scarf's\qss combinatorial\qss theorem\qss in\qss an\qss
abstract\qss setting}{scarf}

\vspace{2.25pt}
\myuppar{Notations.}
Let\sss $A$\sss be set\halfff.\oss
For\qss $a\qff \in\qff A$\qss
we denote by\qss
$A\qff -\qff a$\qss the set\qss 
$A\qff \smallsetminus\qff \{\qff a \qff\}$\nnsp,\oss
and\dss for\qss $b\qff \not\in A$\qss
wel denote by\qss
$A\qff +\qff b$\qss
the set\qss
$A\qff \cup\qff \{\qff b \qff\}$\nnsp.\oss
The set\qss $A\qff -\qff a$\qss is defined only\dss if\qss $a\qff \in\qff A$\nnsp,\oss
and\qss $A\qff +\qff b$\qss is defined only\dss if\qss $b\qff \not\in\qff A$\nnsp.\oss
By\dss $\num{A}$ we denote the number of\dss elements of\dss $A$\nnsp.

\myuppar{Linear orders and\dss dominant sets.}
Let $T$ be a finite set\halfff.\oss
Suppose that a family of\dss linear orders $<_{\dff i}$ on\dss $T$\dnsp,\oss
labeled\dss by elements $i$ of\dss a\sss finite set $I$\nnsp,\pss
is given.\oss
For a non-empty\dss subset\qss $\sigma\qff \subset\qff T$\qss let\qss
$\min_{\dff i}\qff \sigma$\qss 
be the minimal element of $\sigma$ with respect to the order\dss $<_{\dff i}$\nnsp.\oss
A subset\qss $\sigma\qff \subset\qff T$\dss is\dss said\dss to\dss be\qss 
\emph{dominant\dss with\dss respect\dss to}\qss
a non-empty subset\dss 
$C$\dss of\qss $I$\qss if\trs\vspace{4.5pt} 
\begin{equation*}
\quad
\mbox{there\dss is\dss no\dss element}\qff\oss 
y\qff \in\qff T\qff\oss 
\mbox{such\dss that}\qff\oss 
\min\nolimits_{\dff i}\qff \sigma\off <_{\dff i}\off y\qff\oss 
\mbox{for\dss all}\qff\oss
i\qff \in\qff C\dff.
\end{equation*}

\vspace{-7.5pt}
It is convenient\dss to agree that\qss $\varnothing\qff \subset\qff T$\qss
is dominant with respect to every non-empty\qss $C\qff \subset\qff I$\nnsp.\oss
The basic properties of\dss dominant\sss sets are summarized\dss in\dss
the following\dss lemma.\oss

\mypar{Lemma.}{minima}
\emph{If\pss $\sigma\qff \subset\qff T$\dss is\dss 
dominant\dss with\dss respect\dss to\qss $C\qff \subset\off I$\nnsp,\oss
then\qss
$\dis
\sigma
\off =\off
\bigl\{\qff 
\min\nolimits_{\dff i}\dff \sigma \off \bigl|\bigr.\off i\qff \in\qff C 
\qff\bigr\}$\qss
and\dss hence\qss
$\num{\sigma}\off \leq\off \num{C}$\nnsp.\oss
If\oss $\tau\off \subset\off \sigma$\nnsp,\oss
then\dss $\tau$\dss is\dss
dominant\dss with\dss respect\dss to\dss $C$\nnsp,\oss
and\trs if\oss $C\off \subset\off D$\nnsp,\oss
then\dss $\sigma$\dss is\dss
dominant\dss with\dss respect\dss to\dss $D$\nnsp.\oss}

\proof
Clearly\halfff,\pss $\sigma$ contains all minima\dss $\min_{\fff i}\dff \sigma$\nnsp.\oss
Suppose that\qss $x\qff \in\qff \sigma$\qss is different from all\qss
$\min_{\fff i}\dff \sigma$\qss with\qss $i\qff \in\qff C$\nnsp.\oss
Then\qss $\min_{\fff i}\dff \sigma\qff <_{\fff i}\qff x$\qss for all\qss $i\qff \in\qff C$\nnsp,\oss
contrary to the assumption.\oss 
This proves\sss the\sss first\sss statement\sss of\trs the lemma.\oss
The second\sss statement\sss follows from\dss the fact\dss that\qss
$\min\nolimits_{\dff i}\qff \sigma
\off \leq_{\dff i}\off
\min\nolimits_{\dff i}\qff \tau$\qss
for every\dss $i\qff \in\pff I$\nnsp,\oss
and\dss the\sss third\sss statement\dss is\dss 
obvious.\oss \eproof

\myuppar{Cells,\oss rooms,\oss and\sss doors.}
A\dss \emph{cell}\oss is\dss defined as\sss a\sss pair 
$(\dff \sigma\fff,\qff C \trf)$ of\dss subsets\sss
$\sigma\qff \subset\qff T$\sss and\sss 
$C\qff \subset\pff I$\sss
such\dss that\sss $C$ is\dss non-empty\sss  
and
$\sigma$ is\dss dominant with respect\dss to\dss $C$\nnsp.\oss
If\dss $(\dff \sigma\fff,\qff C \trf)$\sss is\dss a cell\halfff,\oss
then\qss 
$\num{\sigma}\qff \leq\qff \num{C}$\qss by\qss Lemma\dss \ref{minima}.\oss  
A cell\sss $(\dff \sigma\fff,\qff C \trf)$\sss 
is\dss called\dss a\qss
\emph{room}\oss if\pss
$\num{C}
\off =\off 
\num{\sigma}$\nnsp,\oss
and\sss a\qss \emph{door}\oss if\pss
$\num{C}
\off =\off 
\num{\sigma}\qff +\qff 1$\nnsp.\oss

A\dss pair\sss $(\dff \tau\fff,\qff D \trf)$\sss 
of\dss subsets\sss
$\tau\qff \subset\qff T$\sss and\sss 
$D\qff \subset\pff I$\sss
is\dss said\sss to\sss be a\qss \emph{door\qss of}\qss a\sss cell\dss 
$(\dff \sigma\fff,\qff C \trf)$\dss
if\qss either\dss 
$(\dff \tau\fff,\qff D \trf)
\off =\off
(\dff \sigma\qff -\qff x\fff,\qff C \trf)$\dss
for\dss some\qss $x\qff \in\qff \sigma$\nnsp,\qff\oss 
or\dss
$(\dff \tau\fff,\qff D \trf)
\off =\off
(\dff \sigma\fff,\qff C\qff +\qff i \trf)$\dss
for some\qss 
$i\qff \in\pff I\qff \smallsetminus\qff C$\nnsp.\oss
By\qss Lemma\qss \ref{minima}\qss every\sss such\dss pair\sss $(\dff \tau\fff,\qff D \trf)$\sss
is\dss a\sss  
a\sss cell\halfff.\oss
It\dss follows\dss that\sss every\dss room\dss has\dss $\num{I}$\dss doors.\oss

Equivalently\halfff,\pss a\sss pair\dss
$(\trf \tau\fff,\qff D \trf)$\dss 
is\trs a\dss door\dss of\dss a\sss cell\dss
$(\trf \sigma\fff,\qff C \trf)$\dss
if\qss 
either\dss 
$(\dff \sigma\fff,\qff C \trf)
\off =\off
(\trf \tau\qff +\qff y\fff,\qff D \trf)$\dss
for some\dss $y\qff \in\pff T\qff \smallsetminus\qff \sigma$\nnsp,\oss 
or\qss 
$(\dff \sigma\fff,\qff C \trf)
\off =\off
(\dff \tau\fff,\qff D\qff -\qff j \trf)$\dss
for some\qss $j\qff \in\pff D$\nnsp.\oss

Clearly\halfff,\oss if\dss
$(\dff \tau\fff,\qff D \trf)$\dss 
is\dss a\sss door\dss of\trs
$(\dff \sigma\fff,\qff C \trf)$\nnsp,\oss
then\dss
$\num{D}\qff -\qff \num{\tau}
\off =\off
\num{C}\qff -\qff \num{\sigma}
\qff +\qff
1$\nnsp.\oss
It\dss follows\sss that\dss if\dss a\sss cell\dss
$(\dff \sigma\fff,\qff C \trf)$\dss
is\dss a\sss room,\oss
then\dss $(\dff \sigma\fff,\qff C \trf)$\dss
is\dss not\sss a\sss door of\dss any\sss cell\halfff.\oss
Similarly\halfff,\oss if\dss $(\dff \tau\fff,\qff D \trf)$\dss 
is\dss a\sss door of\trs a cell\dss
$(\dff \sigma\fff,\qff C \trf)$\nnsp,\oss
then\sss $(\dff \sigma\fff,\qff C \trf)$\sss
is\dss a\sss room\halfff.\oss
A door\dss
$(\trf \tau\fff,\qff D \trf)$\dss
is\dss called\sss an\qss \emph{outside\sss door}\oss
if\qss $\tau\off =\off \varnothing$\qss
and\dss an\qss \emph{internal\dss door}\pss otherwise.\oss

One can\dss imagine\sss that\dss $T$\sss together\dss with\dss the rooms and doors\dss
is\dss some sort\sss of\dss a\sss building\halfff.\oss
As we will\sss see,\oss an outside door\dss is\dss a\sss door
of\dss exactly\sss one room,\oss
and\sss an\dss internal\sss door\dss is\dss a\sss common\sss door\sss of\dss
exactly\dss two rooms.\oss
These properties are\sss the main\sss properties of\trs this\sss building\halfff.\oss

\mypar{Lemma.}{outside-doors}\oss
\emph{Suppose\sss that\qss
$\num{D}\off =\off 1$\nnsp.\oss
Then\dss
$(\trf \varnothing\fff,\qff D \trf)$\dss
is\dss an outside door\dss and\dss
is\dss a door\sss of\dss exactly\sss one room.\oss
Every\sss outside door\dss has\sss this\dss form.\oss}\vspace{-0.75pt}

\proof
Since $\varnothing$ is\dss dominant\dss with\dss respect\dss to every\dss
non-empty\sss subset\sss of\trs $I$\nnsp,\oss
the pair\dss $(\trf \varnothing\fff,\qff D \trf)$\dss is\dss a cell\halfff.\oss
Since\qss $\num{D}\off =\off 1$\nnsp,\oss this cell\dss
is\dss a\sss door\sss and\dss hence an outside door\halfff.\oss
Also,\pss $\num{D}\off =\off 1$\qss implies\sss that\trs
$D\qff =\qff \{\trf i \qff\}$\dss for some\dss $i\qff \in\pff I$\nnsp.\oss
If\qss $(\trf \varnothing\fff,\pff  \{\trf i \qff\} \trf)$\qss
is\dss a\dss door of\dss a cell\dss $(\dff \sigma\fff,\pff C \trf)$\dnsp,\oss
then\qss either\dss
$(\dff \sigma\fff,\pff C \trf)
\off =\off
(\dff \varnothing\fff,\pff  \{\trf i \qff\}\qff +\qff j \trf)$\dss
for some\dss $j\qff \not\in\qff \{\trf i \qff\}$\nnsp,\oss
i.e.\qss for some\dss
$j\off \neq\off i$\nnsp,\oss
or\dss$(\trf \sigma\fff,\pff C \trf)
\off =\off
(\trf \{\trf x \qff\}\fff,\pff  \{\trf i \qff\} \trf)$\dss
for some\qss $x\qff \in\qff T$\nnsp.\oss
In\dss the first\sss case\dss
$\num{\sigma}\off =\off 0$\dss
and\dss
$\num{C}\off =\off 2$\dss
and\dss hence\dss $(\trf \sigma\fff,\pff C \trf)$\dss 
is\dss not\sss a\sss room.\oss
In\dss the second case
$(\trf \sigma\fff,\pff C \trf)$
is\dss a\sss cell\qss if\trs and\dss only\trs if\dss
$\{\trf x \qff\}$ is dominant with respect to\dss $\{\trf i \qff\}$\nnsp,\oss
i.e.\qss if\trs and\dss only\trs if\sss
$x$\sss is\dss the
maximal\dss element\sss of\trs $T$\dss  with\sss respect\dss to\dss $<_{\dff i}$\nsp.\oss
It\dss follows\dss that\sss 
$(\trf \varnothing\fff,\pff  \{\trf i \qff\} \trf)$\sss
is\dss a\dss door\sss of\trs exactly\sss 
one room.\oss
Clearly\halfff,\pss
$(\trf \varnothing\fff,\qff D \trf)$\dss is\dss a\sss door
only\trs if\qss $\num{D}\off =\off 1$\nnsp.\oss  \eproof\vspace{-0.75pt}

\mypar{Lemma.}{doors}\oss
\emph{An\dss internal\sss door\dss $(\trf \tau\fff,\qff D \trf)$\sss
is\dss a\trs door\sss of\dss exactly\dss two\sss rooms.\oss}

\proof
Since\dss $(\trf \tau\fff,\qff D \trf)$\dss
is\dss an\dss internal\sss door\halfff,\oss
$\num{\tau}\qff \geq\qff 1$\dss
and\dss
$\num{D}\off =\off \num{\tau}\qff +\qff 1\off \geq\off 2$\nnsp.\oss
Suppose\sss that\sss $(\trf \tau\fff,\qff D \trf)$ 
is\dss a\sss door of\sss $(\trf \sigma\fff,\qff C \trf)$\dnsp.\oss
Then\sss either\dss
$(\trf \sigma\fff,\qff C \trf)
\off =\off
(\trf \tau\qff +\qff x\fff,\qff D \trf)$\dss
for some\dss $x\qff \not\in\pff \tau$\dnsp,\oss
or\qss
$(\trf \sigma\fff,\qff C \trf)
\off =\off
(\trf \tau\fff,\qff D\qff -\qff i \trf)$\dss
for some\dss $i\qff \in\pff D$\nnsp.\oss
Since\qss
$\num{D}\off \geq\off 2$\nnsp,\oss
the set\sss $C$\sss is\dss
non-empty\halfff.\oss
Therefore,\oss if\dss $\sigma$ is\dss dominant\dss with\dss respect\dss to $C$\nnsp,\oss
then\sss 
$(\trf \sigma\fff,\qff C \trf)$\sss 
is\dss a\sss cell\sss and\dss hence\dss
is\dss a\sss room.\oss
Let\dss us\dss find\sss out\dss when
$\tau\qff +\qff x$\dss is\dss dominant\dss with\dss respect\dss to\dss $D$\sss
and\sss 
$\tau$\dss is\dss dominant\dss with\dss respect\dss to\dss $D\qff -\qff i$\nnsp.
Since\dss 
$\num{\tau}\off =\off \num{D}\qff +\qff 1$\nnsp,\oss
Lemma\qss \ref{minima}\qss implies that there is a unique pair\qss
$\{\dff a\fff,\pff b \qff\}\qff \subset\pff D$\qss
such\dss that\vspace{1.5pt}
\[
\quad
\min\nolimits_{\dff a}\dff \tau
\off =\off
\min\nolimits_{\dff b}\dff \tau
\]

\vspace{-10.5pt}
and\qss $a\qff \neq\qff b$\nnsp.\qff\oss
For\oss $i\off =\off a$\qss or\qss $b$\oss let\oss\vspace{1.5pt}
\[
\quad
\mathbb{M}_{\dff i}
\off\off =\off\off
\bigl\{\off y
\qff \in\qff T 
\off \bigl|\bigr.\off 
\min\nolimits_{\dff k}\dff \tau
\off <_k\off y
\off\mbox{for\dss all}\off
k\qff \in\qff D\qff -\qff i
\off\bigr\}\dff.
\]

\vspace{-10.5pt}
When\pss $\mathbb{M}_{\dff i}\off \neq\off \varnothing$\nnsp,\oss
we will\sss denote\sss by $m_{\dff i}$ the maximal\sss element\sss of\qss $\mathbb{M}_{\dff i}$\qss
with respect\dss to\qss $<_{\dff i}$\nsp.\oss
Since $\tau$ is\dss dominant\dss with\dss respect\dss to\sss $D$\nnsp,\oss
the\sss intersection\dss
$\mathbb{M}_{\dff a}\qff \cap\pff \mathbb{M}_{\dff b}$\dss
is\dss empty\halfff.\oss
Therefore,\oss
if\qss both\dss $\mathbb{M}_{\dff a}$\dss and\dss  $\mathbb{M}_{\dff b}$\dss
are non-empty\halfff,\oss
then\dss
$m_{\dff a}\off \neq\off m_{\dff b}$\dss
and\dss hence\sss the\sss pairs\dss
$(\trf \tau\qff +\qff m_{\dff a}\dff,\pff D \trf)$\dss
and\dss
$(\trf \tau\qff +\qff m_{\dff b}\dff,\pff D \trf)$\dss
are different\halfff.\oss
This reduces\sss the\sss lemma\sss to\sss the following\dss two sublemmas.

\mylemma{Sublemma.}{door-of-minus}
\emph{The set\dss 
$\tau$\dss is\dss dominant with respect\dss to\dss $D\qff -\qff i$\dss
if\dss and\dss only\dss if\pss
$i\qff \in\qff \{\dff a\fff,\pff b \qff\}$\dss 
and\qss $\mathbb{M}_{\dff i}\off =\off \varnothing$\dnsp.}\vspace{-1.75pt}

\subproof\qss 
If\pss $i\off \neq\off a\fff,\pff b$\nnsp,\oss
then\dss the set\qss
$\dis
\bigl\{\pff 
\min\nolimits_{\dff k}\dff \tau 
\off |
\off 
k\qff \in\qff D\qff -\qff i 
\pff
\bigr\}
$\qss
has\qss
$\leq\qff
\num{D}\qff -\qff 2
\off =\off
\num{\tau}\qff -\qff 1$\oss
elements\sss
and\dss hence\dss $\tau$\dss
is\dss not\sss dominant\sss with\sss respect\dss to\dss $D\qff -\qff i$\qss
by\qss Lemma\qss \ref{minima}.\oss 
If\qss $i\off =\off a$\dss or\dss $b$\nnsp,\oss 
then\dss $\tau$\dss is\dss dominant\dss with respect\dss to\qss $D\qff -\qff i$\qss
if\trs and\dss only\trs if\pss 
$\mathbb{M}_{\dff i}\off =\off \varnothing$\nnsp.\oss  \esubproof

\mylemma{Sublemma.}{door-of-plus}
\emph{$\tau\qff +\qff x$\qss is\dss dominant\dss with\sss respect\dss to\qss $D$\qss
if\qss and\dss only\qss if\pss
$x\off =\off m_{\dff i}$\pss
for\dss some\qss 
$i\qff \in\pff \{\dff a\fff,\pff b \qff\}$\dss
such\dss that\oss
$\mathbb{M}_{\dff i}\off \neq\off \varnothing$\nnsp.\oss}\vspace{-0.25pt}

\subproof
To begin with,\oss let us observe that\vspace{3pt}
\begin{equation}
\label{minimum-1}
\quad
\min\nolimits_{\dff i}\pff (\trf \tau\qff +\qff x \trf)
\off =\off 
x
\hspace*{4.2em}\mbox{if}\hspace*{1.5em}
x\off <_{\dff i}\off \min\nolimits_{\dff i}\qff \tau\dff,
\hspace*{1.5em}\mbox{and}\hspace*{1.5em}
\end{equation}

\vspace*{-37.5pt}
\begin{equation}
\label{minimum-2}
\quad
\min\nolimits_{\dff i}\pff (\trf \tau\qff +\qff x \trf)
\off =\off 
\min\nolimits_{\dff i}\qff \tau
\hspace*{1.5em}\mbox{if}\hspace*{1.5em}
\min\nolimits_{\dff i}\qff \tau
\off <_{\dff i}\off 
x\dff.
\end{equation}

\vspace{-9pt}
In particular\halfff,\oss
$\min\nolimits_{\dff i}\qff (\trf \tau\qff +\qff x \trf)
\off =\off
\min\nolimits_{\dff i}\qff \tau$\oss
or\oss
$x$\oss 
for every\oss $i\qff \in\qff D$\dnsp.\oss
Lemma\qss \ref{minima}\qss implies\sss that\qss\vspace{3pt}
\[
\quad
\bigl\{\pff 
\min\nolimits_{\dff i}\qff \tau
\off |\off 
i\qff \in\qff D 
\pff\bigr\}
\off =\off\dff
\tau
\hspace*{1.5em}\mbox{and\qss that}\hspace*{1.5em}
\]

\vspace{-36pt}
\[
\quad
\bigl\{\off 
\min\nolimits_{\fff i}\pff  (\dff \tau\qff +\qff x \dff) 
\off \bigl|\bigr.\off 
i\qff \in\qff D 
\off\bigr\}
\off\off
=\off\off
\tau
\off +\off
x 
\]

\vspace{-9pt}
if\pss
$\tau\qff +\qff x$\qss is\dss dominant\dss with\sss respect\dss to $D$\nnsp.\oss
This may happen only\dss if\qss
$\min\nolimits_{\dff i}\qff (\trf \tau\qff +\qff x \trf)
\off =\off
\min\nolimits_{\dff i}\qff \tau$\oss
for all\dss
$i\qff \in\qff D\qff \smallsetminus\qff \{\qff a\fff,\pff b \qff\}$\qss
and\dss for $i$ equal\dss to one of\trs the elements
of\dss the pair\qss $\{\dff a\fff,\pff b \qff\}$\nnsp,\oss 
and\trs if\qss
$\min\nolimits_{\dff i}\qff (\trf \tau\qff +\qff x \trf)
\off =\off
x$\qss
for $i$ equal\dss to\sss the\sss other element\dss of\qss 
$\{\dff a\fff,\pff b \qff\}$\nnsp.\oss

Therefore,\oss 
if\pss
$\tau\qff +\qff x$\qss is\dss dominant\dss with\sss respect\dss to $D$\nnsp,\oss
we may\sss assume\sss that\vspace{3pt}
\begin{equation}
\label{change-of-minima}
\quad
\min\nolimits_{\dff i}\qff (\trf \tau\qff +\qff x \trf)
\off =\off
\min\nolimits_{\dff i}\qff \tau
\hspace*{1.2em}\mbox{for\dss all}\hspace*{1.2em}
i\qff \in\qff D\qff -\qff a
\hspace*{1em}\mbox{and}\hspace*{1em}
\end{equation}

\vspace{-37.5pt}
\begin{equation}
\label{change-of-minima-1}
\quad
\min\nolimits_{\dff a}\qff (\trf \tau\qff +\qff x \trf)
\off =\off
x\dff.
\end{equation}

\vspace{-9pt}
By\qss (\ref{minimum-1})\qss and\qss (\ref{minimum-2})\qss in\dss this case\oss
$\min\nolimits_{\dff i}\qff \tau
\off <_{\dff i}\off 
x$\oss
for all\qss
$i\qff \in\qff D\qff -\qff a$\qss
and\qss
$x
\off <_{\dff a}\off
\min\nolimits_{\dff a}\qff \tau$\dnsp.\oss
It\dss follows\dss that\qss
$x\qff \in\qff \mathbb{M}_{\dff a}$\qss
and\dss that\qss $\tau\qff +\qff x$\qss can\sss be dominant with\sss respect\dss to\dss $D$\dss
only\qss if\dss $x$\dss is\dss the maximal\sss element\sss of\dss
$\mathbb{M}_{\dff a}$\dss
with respect\dss to\qss $<_{\dff a}$\nsp,\oss
i.e.\qss only\qss if\pss
$x\off =\off m_{\dff a}$\nnsp.\oss

Conversely\halfff,\oss
if\halfff,\oss say\halfff,\oss
$\mathbb{M}_{\dff a}\qff \neq\qff \varnothing$\qss
and\pss
$x\qff \in\qff \mathbb{M}_{\dff a}$\nnsp,\oss
then\vspace{2.5pt}
\[
\quad
\min\nolimits_{\dff i}\qff \tau
\off <_{\dff i}\off 
x
\hspace*{1.2em}\mbox{for\dss all}\hspace*{1.2em}
i\qff \in\pff D\qff -\qff a\dff.
\]

\vspace{-9.5pt}
If\trs also\qss
$\min\nolimits_{\dff a}\qff \tau
\off <_{\dff a}\off 
x$\nnsp,\oss
then\dss $\tau$\dss is\dss not\sss dominant\dss with\sss respect\dss to\dss $D$\nnsp,\oss
contrary\dss to\sss the assumption.\oss
Therefore\qss
$x
\off <_{\dff a}\off
\min\nolimits_{\fff a}\qff \tau$\dnsp.\oss
By applying\qss 
(\ref{minimum-1})\qss and\qss (\ref{minimum-2})\qss
again,\oss 
we see that\qss (\ref{change-of-minima})\qss and\qss (\ref{change-of-minima-1})\qss hold.\oss
It\dss follows\dss that\qss
if\pss
$x\off =\off m_{\dff a}$\nsp,\oss
then\dss
$\tau\qff +\qff x$\qss is\dss dominant\dss 
with\dss respect\dss to\sss $D$\nnsp.\oss  \esubproof  \eproof

\myuppar{Colorings.}
A\dss \emph{coloring}\qss is\dss arbitrary\dss map\qss
$T\ttoo I$\nnsp.\oss
Let\dss us\dss fix\sss a\sss coloring\sss $c$\nnsp.\oss
A\sss cell\dss $(\dff \sigma\fff,\qff C \trf)$\dss
is\dss called\qss \emph{colorful}\oss if\oss
$C\off =\off c\dff(\trf \sigma \trf)$\nnsp.\oss
A colorful\sss cell\trs is\dss automatically\dss a\dss room.\oss
Scarf's\qss combinatorial\dss theorem\sss asserts\sss that\dss
for every\sss coloring\dss there exists a\sss 
colorful\dss room.\oss
See\qss Theorem\qss \ref{main-lemma}.\oss

A cell\dss $(\dff \sigma\fff,\qff C \trf)$\dss
is\dss called\qss \emph{nearly\dss colorful}\oss if\oss
$\num{C\qff \smallsetminus\qff c\dff(\trf \sigma \trf)}\off =\off 1$\nnsp.\oss
If\dss $(\dff \sigma\fff,\qff C \trf)$\dss 
is\dss nearly\sss colorful,\oss
then\dss
$C\qff \smallsetminus\qff c\dff(\dff \sigma \trf)$\dss
consists of\dss one element\halfff,\oss
called\trs the\qss \emph{type}\pss of\qss
$(\trf \sigma\fff,\qff C \trf)$\nnsp.\oss

\mypar{Lemma.}{end-types}
\emph{For every\dss $i\qff \in\pff I$\dss 
there\dss is\dss exactly\sss one outside door
of\qss the\sss type\sss $i$\nnsp.\qff\oss
Every\sss door of\dss a colorful\dss room\dss is\dss nearly\sss colorful\dss
and\sss there\dss is\dss exactly\sss 
one door of\qss each\sss type among\dss them.\oss}\vspace{-0.6pt}

\proof
The outside doors are\sss the cells of\trs the form\dss
$(\trf \varnothing\fff,\pff  \{\trf i \qff\} \trf)$\nnsp,\oss
where\qss $i\qff \in\pff I$\nnsp.\oss
Clearly\halfff,\oss the door\dss $(\trf \varnothing\fff,\pff  \{\trf i \qff\} \trf)$\dss
is\dss nearly\sss colorful\sss and\dss its\sss type\dss is\sss $i$\nnsp.\oss
This proves\sss the first\sss statement\halfff.\oss

Let\dss $(\dff \sigma\fff,\qff C \trf)$\dss be a colorful\dss room.\oss
If\qss  
$x\qff \in\qff \sigma$\nnsp,\oss 
then\dss 
$\num{C\qff \smallsetminus\qff c\trf(\dff \sigma\qff -\qff x \trf)}
\qff =\qff
1$\dss
and\dss hence\sss
the door\dss $(\trf \sigma\qff -\qff x\fff,\pff C \trf)$\qss is\dss 
nearly\sss colorful\halfff.\oss
If\qss
$i\qff \in\pff I\qff \smallsetminus\qff C$\nnsp,\oss 
then\dss 
$\num{(\qff C\qff +\qff i\qff)\qff \smallsetminus\qff c\trf(\dff \sigma \dff)}
\qff =\qff
1$\dss
and\dss hence\sss 
the door\dss $(\trf \sigma\fff,\pff C\qff +\qff x \trf)$\qss is\dss 
nearly\sss colorful\halfff.\oss
Obviously\halfff,\oss the\sss type of\qss 
$(\trf \sigma\qff -\qff x\fff,\pff C \trf)$\qss 
is\dss 
$c\dff(\dff x\trf)\qff \in\qff c\dff(\dff \sigma\trf)$\nnsp,\oss
and\dss the\sss type of\qss $(\trf \sigma\fff,\pff C\qff +\qff i \trf)$\qss is 
$i\qff \in\pff I\qff \smallsetminus\qff C$\nnsp.\oss
Since\sss $(\dff \sigma\fff,\qff C \trf)$\sss is\dss colorful\halfff,\oss
the map\sss $c$
induces a\sss bijection\qss
$\sigma\qff \ttoo\qff 
C$\nnsp.\oss
The\sss second\sss statement\dss follows.\oss  \eproof

\mypar{Lemma.}{types}
\emph{Suppose\sss that\dss $(\trf \tau\fff,\qff D \trf)$\dss
is\dss a\dss nearly\dss colorful\dss door\dss
of\qss a\sss room\sss $(\dff \sigma\fff,\qff C \trf)$\nnsp.\oss
Then\dss $(\dff \sigma\fff,\qff C \trf)$\dss
is\dss either\sss colorful\halfff,\oss
or\dss nearly\dss colorful\dss of\qss 
the\sss same\dss type\sss as\dss $(\trf \tau\fff,\qff D \trf)$\nnsp.\oss}\vspace{-0.6pt}

\proof
Clearly\halfff,\pss
$C\qff \smallsetminus\qff c\dff(\dff \sigma \trf)
\off \subset\off
D\qff \smallsetminus\qff c\dff(\trf \tau \trf)$\nnsp.\oss
Together\dss with\dss
$\num{D\qff \smallsetminus\qff c\dff(\trf \tau \trf)}
\off =\off 
1$\dss
this\sss implies\sss that\qss
$\num{C\qff \smallsetminus\qff c\dff(\dff \sigma \trf)}
\off \leq\off
1$\nnsp.\oss
Since\dss
$\num{c\dff(\dff \sigma \trf)}
\off \leq\off
\num{\sigma}
\off \leq\off
\num{C}$\nnsp,\qff\oss
if\pss
$\num{C\qff \smallsetminus\qff c\dff(\dff \sigma \trf)}
\off =\off
0$\nnsp,\oss
then\dss $C\off =\off c\dff(\dff \sigma \trf)$\dss
and\dss hence $(\dff \sigma\fff,\qff C \trf)$\sss is\dss colorful\halfff.\oss
If\qss
$\num{C\qff \smallsetminus\qff c\dff(\dff \sigma \trf)}
\off =\off
1$\nnsp,\oss
then\dss
$C\qff \smallsetminus\qff c\dff(\dff \sigma \trf)
\off =\off
D\qff \smallsetminus\qff c\dff(\trf \tau \trf)$\dss 
and\dss hence\sss
$(\dff \sigma\fff,\qff C \trf)$\sss is\dss nearly\sss colorful\sss
and\dss has\sss the same\sss type as $(\trf \tau\fff,\qff D \trf)$\nnsp.\oss  \eproof

\mypar{Lemma.}{basic-cases}
\emph{If\pss $(\dff \sigma\fff,\qff C \trf)$ is\dss a\dss nearly\qss colorful\qss room,\oss
then\dss
$\num{c\dff(\trf \sigma \trf)}$\dss is\trs equal\dss either\dss to\qss $\num{\sigma}$\dss
or\qss to\qss $\num{\sigma}\qff -\qff 1$\nnsp.\oss
If\oss $(\dff \sigma\fff,\qff C \trf)$\dss is\trs a\dss nearly\qss colorful\qss door\halfff,\oss
then\qss
$c\dff(\trf \sigma \trf)\off \subset\off C$\nnsp.\oss}\vspace{-0.6pt}

\proof
If\qss  $(\dff \sigma\fff,\qff C \trf)$ is\dss a\sss nearly\sss colorful\dss room,\oss
then\dss $\num{C}\off \leq\off \num{c\dff(\dff \sigma\dff)}\qff +\qff 1$\dss
and\dss hence\vspace{1.5pt}
\[
\quad
\num{\sigma}
\off \leq\off 
\num{C}
\off \leq\off
\num{c\dff(\dff \sigma\dff)}\qff +\qff 1
\off \leq\off 
\num{\sigma}\qff +\qff 1
\qff.
\]

\vspace{-10.5pt}
It\dss follows\dss that\sss
$\num{c\dff(\dff \sigma\dff)}
\off =\off 
\num{\sigma}$\dss or\dss 
$\num{\sigma}\qff -\qff 1$\nnsp.\oss
If\qss  $(\dff \sigma\fff,\qff C \trf)$ is\dss a\sss nearly\sss colorful\dss door\halfff,\oss
then\dss\vspace{1.5pt} 
\[
\quad\num{\sigma}
\off \geq\off
\num{c\dff(\dff \sigma \dff)}
\off \geq\off
\num{c\dff(\dff \sigma \dff)\qff \cap\qff C}
\off =\off
\num{C}\qff -\qff 1
\qff =\qff 
\num{\sigma}
\]

\vspace{-10.5pt}
and\dss hence\dss 
$\num{c\dff(\dff \sigma \dff)}
\off =\off 
\num{c\dff(\dff \sigma \dff)\qff \cap\qff C}$\dss
and\qss
$c\dff(\dff \sigma \dff)
\off =\off 
c\dff(\dff \sigma \dff)\qff \cap\qff C$\nnsp.\oss 
It\dss follows\dss that\qss
$c\dff(\dff \sigma \dff)\qff \subset\qff C$\nnsp.\oss  \eproof

\mypar{Lemma.}{rooms}
\emph{If\pss $(\dff \sigma\fff,\qff C \trf)$\sss is\dss a\sss
nearly\qss colorful\qss room,\oss
then\dss $(\dff \sigma\fff,\qff C \trf)$\sss
has\dss two\dss nearly\qss colorful\qss doors.}

\proof
Lemma\qss \ref{basic-cases}\pss implies\sss that\sss 
either\dss
$\num{c\dff(\trf \sigma \trf)}
\off =\off 
\num{\sigma}$\dss
or\dss
$\num{c\dff(\trf \sigma \trf)}
\off =\off 
\num{\sigma}\qff -\qff 1
\off =\off
\num{C}\qff -\qff 1$\nnsp.\oss

Suppose\sss first\dss that\dss 
$\num{c\dff(\trf \sigma \trf)}
\off =\off 
\num{\sigma}\qff -\qff 1
\off =\off 
\num{C}\qff -\qff 1$\nnsp.\oss
Then\dss
$c\dff(\trf \sigma \trf)\off \subset\off C$\nnsp.\oss
Therefore,\oss
if\dss also\dss $i\qff \not\in\pff C$\nnsp,\oss
then\dss
$\num{(\trf C\qff +\qff i \trf)\qff \smallsetminus\qff c\dff(\trf \sigma \trf)}
\off =\off
2$\dss
and\dss hence\dss
$(\trf \sigma\fff,\qff C\qff +\qff i \trf)$\dss
is\dss not\dss nearly\sss colorful\halfff.\oss
Since
$\num{c\dff(\trf \sigma \trf)}\off =\off \num{\sigma}\qff -\qff 1$\nnsp,\oss
there are different\sss elements\dss $x\fff,\pff y\qff \in\qff \sigma$\dss
such\dss that\dss
$c\dff(\dff x\trf)\off =\off c\dff(\dff y\trf)$\nnsp.\oss
More\-over\halfff,\oss the coloring\dss $c$\dss is\dss injective on\dss both\dss $\sigma\qff -\qff x$\dss
and\dss $\sigma\qff -\qff y$\nnsp.\oss
If\qss $z\qff \in\pff C$\dss and\dss $z\off \neq\off x\fff,\pff y$\nnsp,\oss
then\dss $c\dff(\dff \sigma\qff -\qff z\trf)$\dss is\dss properly\sss
contained\dss in\dss $c\dff(\dff \sigma\trf)$\nnsp.\oss
Since\dss $c\dff(\dff \sigma\trf)$\dss is\dss properly\sss
contained\dss in\dss $C$\nnsp,\oss
it\dss follows\dss that\qss
$\num{C\qff \smallsetminus\qff c\dff(\trf \sigma\qff -\qff z\trf)}
\off =\off
2$\qss
and\dss hence\dss
$(\dff \sigma\qff -\qff z\dff,\qff C \trf)$\dss
is\dss not\dss nearly\sss colorful\halfff.\oss
On\dss the other\dss hand,\oss
if\qss $z\off =\off x$\dss or\dss $y$\nnsp,\oss 
then\dss
$c\dff(\trf \sigma\qff -\qff z\trf)
\off =\off
c\dff(\trf \sigma\trf)$\dss
and\dss hence\dss
$\num{C\qff \smallsetminus\qff c\dff(\trf \sigma\qff -\qff z\trf)}
\off =\off
1$\nnsp.\oss
In\dss this case\dss
$(\dff \sigma\qff -\qff z\dff,\qff C \trf)$\dss
is\dss nearly\sss colorful\halfff.\oss
Therefore\dss
$(\dff \sigma\fff,\qff C \trf)$\dss
has\sss two\sss nearly\sss colorful\sss doors.\oss

Suppose\dss that\dss
$\num{c\dff(\trf \sigma \trf)}\off =\off \num{\sigma}$\nnsp.\oss
Then  
$c$ is\dss injective on $\sigma$\nnsp.\oss
Since $(\trf \sigma\fff,\qff C \trf)$ is\dss nearly\sss col\-or\-ful\halfff,\oss
this\sss implies\sss that\sss 
$\num{c\dff(\dff \sigma \dff)\qff \smallsetminus\qff C}
\off =\off 
1$
and\dss hence\sss there\dss is\dss a\sss unique\sss 
$y\qff \in\qff \sigma$\sss such\dss that\sss 
$c\dff(\dff y\trf)\qff \not\in\qff C$\nnsp.\oss
If\trs $x\qff \in\qff \sigma$\nnsp,\oss
then\dss
$c\dff(\dff x\trf)\qff \not\in\qff c\dff(\trf \sigma\qff -\qff x\trf)$\dss
because $c$ is\dss injective on $\sigma$\nnsp.\oss
If\dss also $c\dff(\dff x\trf)\qff \in\pff C$\nnsp,\oss
then\vspace{3pt}
\[
\quad
C\qff \smallsetminus\qff c\dff(\trf \sigma\qff -\qff x\trf)
\off =\off
(\qff
C\qff \smallsetminus\qff c\dff(\trf \sigma \trf)
\qff)
\qff +\qff 
c\dff(\dff x\trf)
\]

\vspace{-9pt}
and\dss hence\qss
$\num{C\qff \smallsetminus\qff c\dff(\trf \sigma\qff -\qff x\trf)}
\off =\off 
2$\nnsp.\oss
Therefore\sss in\dss this case\dss
$(\dff \sigma\qff -\qff x\fff,\pff C \dff)$\dss
is\dss not\dss nearly\sss colorful\halfff.\oss

On\dss the other\dss hand,\oss
if\dss $c\dff(\dff x\trf)\qff \not\in\pff C$\nnsp,\oss
i.e.\qss if\pss $x\off =\off y$\nnsp,\oss
then\vspace{4.5pt}
\[
\quad
C\qff \smallsetminus\qff c\dff(\trf \sigma\qff -\qff x\trf)
\off =\off
C\qff \smallsetminus\qff c\dff(\trf \sigma \trf)
\]

\vspace{-7.5pt}
and\dss hence\dss
$\num{C\qff \smallsetminus\qff c\dff(\trf \sigma\qff -\qff x\trf)}
\off =\off
\num{C\qff \smallsetminus\qff c\dff(\trf \sigma \trf)}
\off =\off 
1$\nnsp.\oss
Therefore\dss
$(\dff \sigma\qff -\qff x\fff,\pff C \dff)$
is\dss nearly\sss colorful\qss if\trs and\dss only\trs if\dss
$x\off =\off y$\nnsp.\qff\oss
If\pss $i\qff \in\qff I\qff \smallsetminus\qff C$\qss and\pss
$i\qff \not\in\qff c\dff(\trf \sigma \trf)$\dnsp,\oss
then\vspace{4.25pt}
\[
\quad
(\qff
C\qff +\qff i
\qff)
\qff \smallsetminus\qff c\dff(\trf \sigma \trf)
\off =\off
(\qff
C\qff \smallsetminus\qff c\dff(\trf \sigma \trf)
\qff)
\qff +\qff i
\]

\vspace{-7.5pt}
and\dss hence\qss
$\num{(\dff
C\qff +\qff i
\trf)
\qff \smallsetminus\qff c\dff(\trf \sigma \trf)}
\off =\off
\num{C\qff \smallsetminus\qff c\dff(\trf \sigma \trf)}
\qff +\qff 
1
\qff =\qff 
2$\nnsp.\oss
Therefore\sss in\dss this case\dss
$(\dff \sigma\fff,\pff C\qff +\qff i \trf)$\dss
is\dss not\dss nearly\sss colorful\halfff.\oss
On\dss the other\dss hand,\oss
if\qss $i\qff \in\qff c\dff(\trf \sigma \trf)$\nnsp,\oss
i.e.\pss if\pss $i\off =\off c\trf(\trf y \trf)$\nnsp,\oss
then\vspace{4.5pt}
\[
\quad
(\qff
C\qff +\qff i
\qff)
\qff \smallsetminus\qff c\dff(\trf \sigma \trf)
\off =\off
C\qff \smallsetminus\qff c\dff(\trf \sigma \trf)
\]

\vspace{-7.5pt}
and\dss hence\qss
$\num{(\qff
C\qff +\qff i
\qff)
\qff \smallsetminus\qff c\dff(\trf \sigma \trf)}
\off =\off
\num{C\qff \smallsetminus\qff c\dff(\trf \sigma \trf)}
\off =\off 
1$\nnsp.\oss
Therefore\dss
$(\dff \sigma\fff,\pff C\qff +\qff i \trf)$
is\dss nearly\sss colorful\qss if\trs and\dss only\trs if\dss
$i\off =\off c\trf(\trf y \trf)$\nnsp.\qff\oss
Hence\sss in\dss this case
$(\dff \sigma\fff,\qff C \trf)$\dss
also\sss has\sss two\sss nearly\sss colorful\sss doors.\oss  \eproof

\mypar{Scarf's\qss combinatorial\trs theorem.}{main-lemma}
\emph{There\dss is\dss at\dss least\dss one\sss 
colorful\dss room.\oss
More\-o\-ver\halfff,\oss the number of\dss such\dss rooms\dss is\dss odd.\oss}

\prooftitle{An\dss informal\dss proof\halfff}
The following\dss informal\sss explanation of\dss the proof\dss is\dss based\dss
on\qss Scarf's\qss outline\qss (see\qss \cite{sc3},\oss p.\qss 48).\oss
Actually\halfff,\oss the\sss terms\qss \emph{rooms}\qss and\qss \emph{doors}\qss
were suggested\dss by\dss this outline,\oss
although\qss Scarf\qss used\dss them only\dss informally\halfff.\oss
Let\dss us\dss think\sss about $T$ and\dss its\sss rooms and doors\dss
as\sss if\qss it\dss is\dss a\sss building\halfff.\oss
The above\sss lemmas\sss
tell\sss us all\dss what\dss we need\dss to know about\dss this building\halfff.\oss
We are interested only\dss in colorful\sss and\dss nearly\sss colorful\dss rooms
and\dss in\dss nearly\sss colorful\sss doors.\oss
One can\dss imagine\sss that\sss
all\sss other doors and\dss rooms are permanently\dss locked.\oss

Our\dss plan\dss is\dss to enter\dss this building\dss through an outside door
and\dss then follow rooms and doors,\oss
never\dss turning\dss back.\oss
An outside door\dss leads\sss to a\sss particular\sss room,\oss
which\dss is\dss either colorful\sss or\sss 
has only\sss one other\dss unlocked\sss door\halfff.\oss
This unlocked door\dss leads us\sss to another\sss room,\oss
which\dss is\dss also either colorful\sss or\sss 
has only\sss one other\dss unlocked\sss door\halfff.\oss
We will\sss continue\sss to explore\sss the building\dss in\dss this way\dss until\dss 
we enter a\sss colorful\dss room.\oss
Since we do not\dss turn\dss back,\oss we cannot\dss pass\sss through\dss
the same door or enter\dss the same room\dss twice\qss
(a room\dss has only\dss two unlocked doors).\oss
Potentially\halfff,\oss our journey\dss may\sss end at\sss another outside door\halfff.\oss
But\qss Lemma\qss \ref{types}\qss implies\sss that\dss the\sss types of\dss doors
and\dss rooms stay\dss the same during our\sss journey\halfff,\oss
and\dss by\qss Lemma\qss \ref{end-types}\qss
there\dss is\dss only\sss one outside door of\dss any\sss given\dss type.\oss
Since\sss we cannot\dss return\dss to\sss this door and\dss 
there are only\dss finitely\dss many\dss rooms,\oss
we will\sss eventually\dss reach a colorful\dss room.\oss

We found a particular colorful\dss room.\oss
If\dss there are other\halfff,\oss
then\dss the same argument\sss shows\sss that\dss 
they\sss are pairwise related\dss by\dss routes consisting of\dss
nearly\sss colorful\dss rooms and doors of\trs the\sss type $i$\nnsp,\oss
where $i$ is\dss the\sss type of\trs the outside door used\dss to
enter\dss the building\halfff.\oss
Therefore\sss the number of\dss other colorful\dss rooms\dss is\dss even
and\dss the\sss total\dss number of\dss colorful\dss rooms\dss is\dss odd.\oss

\myuppar{The graphs of\dss rooms and doors.}
In order\dss to\sss turn\dss the above informal\dss proof\dss
in a formal\sss one,\oss it\dss is\dss convenient\dss to represent\dss
journeys\sss through\dss the building\dss by\dss paths in appropriate graphs.\oss
Let\dss us\dss fix\sss some element\qss $i\qff \in\pff I$\qss
and\sss restrict\sss our\sss attention\dss by\dss the\sss
colorful\dss rooms\sss and\dss by\dss the 
nearly\sss colorful\dss rooms and doors of\dss the\sss type $i$\nnsp.\oss
The relation\dss\vspace{3pt}
\[
\quad
(\trf \tau\fff,\qff D \trf)
\hspace*{0.8em}\mbox{is\dss a\dss door\sss of}\hspace*{0.8em}
(\trf \sigma\fff,\qff C \trf)
\]

\vspace{-9pt}
between such\dss rooms and doors can\sss be encoded\dss
in\dss terms of\dss a\sss graph\dss $G_{\dff i}$\nsp,\oss
having as\sss the vertices\sss the\sss colorful\dss rooms and\dss the 
nearly\sss colorful\dss rooms and doors of\dss the\sss type $i$\nnsp.\oss
Each\sss such\sss room\dss is\dss connected\dss by\sss an edge\sss to each
of\qss its nearly\sss colorful\sss doors.\oss 
There are no other edges.\oss
Clearly\halfff,\oss Scarf's\qss combinatorial\dss theorem\dss 
immediately\dss follows\dss from\dss the following\dss fact\halfff.\oss

\mypar{Theorem.}{graph}
\emph{The\dss graph\dss $G_{\dff i}$\sss
consists of\qss several\dss disjoint\trs paths and cycles.\oss
With\sss one exception,\oss
the endpoints\sss of\qss these paths\dss are\sss properly\sss colored\sss cells.\oss
The only\sss exception\dss is\dss
$(\dff \varnothing\fff,\pff  \{\trf i \qff\} \trf)$\nnsp.\oss}

\proof
By\qss Lemma\qss \ref{end-types}\qss
there\dss is\dss only\sss one outside door\dss with\dss the\sss type $i$\nnsp.\oss
By\qss Lemma\qss \ref{outside-doors}\qss it\dss is\dss a door\dss
of\dss exactly\sss one room,\oss
which\dss is\dss nearly\sss colorful\sss and\dss has\sss the\sss type $i$\sss by\qss Lemma\qss \ref{types}.\oss
Therefore\sss this outside door\dss
is\dss an\sss endpoint\sss of\dss exactly\sss one edge.\oss
Also\sss by\qss Lemma\qss \ref{end-types}\qss 
every\sss colorful\dss room\dss 
has exactly\sss one nearly\sss colorful\sss door with\dss the\sss type $i$ and\dss hence\dss
is\dss
also an endpoint\sss of\dss exactly\sss one edge.\oss
Lemmas\qss \ref{doors}\qss and\qss \ref{types}\qss
imply\dss that\sss
every\dss internal\sss door\dss 
is\dss an\dss endpoint\dss of\trs two edges.\oss
Finally\halfff,\oss Lemmas\qss \ref{rooms}\qss and\qss \ref{types}\qss imply\dss that\sss
every\dss nearly\sss colorful\dss room\dss
is\dss an endpoint\sss of\trs two edges.\oss
To sum\dss up,\oss every\dss vertex of\trs $G_{\dff i}$\sss
is\dss an endpoint\sss of\dss one or\dss two edges.\oss
A vertex\dss is\dss an endpoint\sss of\dss only\sss one edge\trs
if\qss and\dss only\trs if\trs it\dss is\dss either\dss
a\sss colorful\dss room\sss or\dss the unique outside door\dss
of\trs the\sss type $i$\nnsp.\oss 
The\sss theorem\dss follows.\oss  \eproof

\myuppar{Remarks.}
The above proof\dss can\sss be easily\dss turned\dss into an algorithm for finding
a\sss properly\sss colored\dss room.\oss
In\dss this circle of\dss questions such\qss 
\emph{path-following\dss algorithms}\qss originate\sss in\dss the work of\qss
Scarf\qss \cite{sc2}\qss and are very\dss popular\halfff.\oss
Such algorithms were actually\dss used,\oss already\dss by\qss Scarf\halfff,\oss
to compute approximations\sss to fixed\dss points 
of\dss self-maps of\dss a simplex\halfff.\oss
See\qss \cite{sc2},\oss \cite{sc3}.\oss

\mysection{Brouwer's\qss fixed\qss point\qss theorem}{brouwer}

\myuppar{The standard simplex\halfff.}
Let\dss us\sss turn\dss to\sss the\sss 
non-combinatorial\dss part\sss of\dss the proof\dss of\qss
Brouwer's\qss theorem.\oss
Let\dss us\dss fix\sss a non-negative integer $n$\nnsp,\pss
and\dss let\qss
$I
\off =\off 
\{\qff 0\fff,\pff 1\fff,\pff \ldots\fff,\pff n \qff\}$\nnsp.\oss
Let\dss us\dss number\dss the coordinates\sss in\dss $\rrr^{\fff n\dff +\dff 1}$\dss
by\sss elements of\trs $I$\nnsp.\oss
For\sss a\sss point\qss 
$x\qff \in\qff \rrr^{\fff n\dff +\dff 1}$\qss
let\dss $x_{\fff i}$\dss be\sss the $i$\dnsp-th coordinate of\dss $x$\nnsp,\oss
so\sss that\oss
$x\off =\off
(\qff x_{\trf 0}\fff,\pff  x_{\dff 1}\fff,\pff \ldots\fff,\pff  x_{\dff n} \qff)$\dnsp.\oss
Let\dss $\Delta^{n}\off \subset\off \rrr^{\fff n\dff +\dff 1}$\dss
be\sss the standard $n$\dnsp-simplex\sss defined\dss by\dss the equation\dss  
$x_{\trf 0}\qff +\qff x_{\dff 1}\qff +\qff \ldots\qff +\qff x_{\dff n}
\off =\off 
1$\dss
and\dss the inequalities\dss $x_{\dff i}\qff \geq\qff 0$\dss
with\dss $i\qff \in\qff I$\nnsp.\oss
Let\dss $l\qff \geq\qff 1$ be\sss another\sss integer\sss  
and\dss let\dss 
$T\off =\off T_{\fff l}$\dss be\sss the set\sss of\dss all\dss 
$x
\qff \in\qff
\Delta^{n}$\qss
such that\dss every\dss $x_{\fff i}$\dss 
is\dss an\dss integer\dss multiple\sss of\dss $1/\fff l$\nnsp.\oss
The set\dss $T$\dss will\sss serve as a discrete approximation\dss to\dss $\Delta^n$\dnsp.\oss

\myuppar{The\sss linear orders on\dss $T$\dnsp.}
For each\qss
$i\qff \in\qff I$\oss
let us choose a linear order\dss $<_{\dff i}$\dss on\dss $T$\dss 
such\dss that\vspace{2.25pt}
\begin{equation}
\label{coordinate-orders}
\quad
x_{\fff i}\off <\off y_{\fff i}
\hspace*{1.2em}
\mbox{implies}
\hspace*{1.2em}
x\off <_{\fff i}\off y
\end{equation}

\vspace{-9.75pt}
for every\qss $x\fff,\pff y\qff \in\pff T$\qss
(obviously\halfff,\pss such orders exist\fff).\oss

\mypar{Theorem.}{dominant-size}
\emph{Let\pss 
$\sigma\qff \subset\pff T$\qss and\pss $C\qff \subset\pff I$\nnsp.\oss
If\oss
$\sigma$\dss
is\dss dominant\dss with\dss respect\dss to\qss $C$\nnsp,\oss 
then}\vspace{2.25pt}
\[
\quad
|\qff x_{\dff i}\qff -\qff y_{\dff i} \qff|
\off <\off
2\trf (\dff n\qff +\qff 1\dff)\left/\dff l\right.
\]

\vspace{-9.75pt}
\emph{for\dss every\oss $x\fff,\pff y\qff \in\qff \sigma$\qss 
and\qss $i\qff \in\pff I$\qss
and\qss
$x_{\dff i}\off <\off (\dff n\qff +\qff 1\dff)\fff/\fff l$\qss
for\dss every\oss $x\qff \in\qff \sigma$\qss 
and\qss $i\qff \in\qff I\qff \smallsetminus\qff C$\nnsp.\oss}

\proof
For each\qss $i\qff \in\pff I$\qss
let\qss
$m\dff(\dff i\trf)
\off =\off 
\min\nolimits_{\dff i}\qff \sigma$\hnsp.\oss
Let\qss
$m_{\dff i}
\off =\off 
m\dff(\dff i\trf)_{\dff i}$\qss
be\sss the $i${\nnsp}th\dss coordinate of\dss $m\dff(\dff i\trf)$\dss
for\qss $i\qff \in\pff C$\qss
and\qss let\qss
$m_{\dff i}
\off =\off
0$\qss
for\qss
$i\qff \in\qff I\qff \smallsetminus\qff 
C$\nnsp.\oss
By\dss using\dss the\sss triangle inequality\sss and\dss the fact\dss that\qss
$\num{C}\qff \leq\qff n\qff +\qff 1$\nnsp,\oss
we see\sss 
that\dss it\dss is\dss sufficient\dss to
prove\sss that\vspace{1.5pt}
\begin{equation*}
\quad
0
\off \leq\off 
x_{\dff i}\qff -\qff m_{\dff i}
\off <\off 
\num{C}\fff\left/\dff l\right.
\end{equation*}

\vspace{-10.5pt}
for\dss every\qss $x\qff \in\qff \sigma$\qss and\qss $i\qff \in\pff I$\nnsp.\oss
The\sss inequalities\qss
$0
\qff \leq\qff 
x_{\dff i}\qff -\qff m_{\dff i}$\qss
hold\dss by\dss 
the definition of\dss $m_{\dff i}$\nnsp.\oss
As\sss the first\sss  step\sss toward\dss the 
inequalities\qss
$x_{\dff i}\qff -\qff m_{\dff i}
\off <\off 
\num{C}\fff\left/\dff l\right.$\nnsp,\qff\oss
let\dss us\dss prove\sss that\vspace{3pt}
\begin{equation}
\label{c-sum}
\quad 
1\off -\off
\sum_{k\qff \in\qff C\vphantom{K^K}}\qff m_{\dff k}
\off\qff <\off\qff
\num{C}\fff\left/\dff l\right. 
\qff.
\end{equation}

\vspace{-12pt}\vspace{0.25pt}
If\trs this\dss is\dss not\dss the case,\oss then\off\oss\vspace{-1.5pt}
\[
\quad
\sum_{k\qff \in\qff C\vphantom{K^K}}\qff
\bigl(\qff 
m_{\dff k}
\qff +\qff 
(\dff 1\fff/\fff l \qff)
\qff\bigr)
\off\qff =\off\off
\left(\off 
\sum_{k\qff \in\qff C\vphantom{K^K}}\qff m_{\dff k}
\off\right)
\off +\off\dff 
\num{C}\fff\left/\dff l\right.
\off\off \leq\off\off
 1
\qff.
\]

\vspace{-10.5pt}
Since every\sss $m_{\fff i}$\sss is\dss a\sss multiple of\trs $1/\fff l$\nnsp,\oss
This\sss implies\sss that\dss there exists 
a\sss point\qss $M\qff \in\qff T$\qss such that\qss
$M_{\dff k}
\off\off \geq\off\off
m_{\dff k}\qff +\qff (\dff 1/\fff l \dff)$\qss
and\dss hence\dss
$\min\nolimits_{\dff k}\dff \sigma
\off <_{\dff k}\off
M$\dss
for every\dss $k\qff \in\qff C$\nnsp.\oss
The contradiction\dss with $\sigma$
being\sss dominant\dss with\sss 
respect\dss to\dss $C$\dss
proves\sss (\ref{c-sum}).\oss

Let\qss
$x\qff \in\qff \sigma$\qss and\qss $i\qff \in\pff I$\nnsp.\qff\oss
If\qss $i\qff \in\pff C$\nnsp,\qff\oss then\vspace{3pt} 
\[
\quad
x_{\dff i}\qff -\qff m_{\dff i}
\off\off \leq\off\off 
\sum_{k\qff \in\qff C\vphantom{K^K}}\qff 
\left(\qff 
x_{\dff k}\qff -\qff m_{\dff k}
\qff\right)
\off\off \leq\off\off
\sum_{k\qff \in\pff I\vphantom{K^K}}\qff x_{\dff k}
\off -\off 
\sum_{k\qff \in\qff C\vphantom{K^K}}\qff m_{\dff k}
\off\off =\off\off 
1
\off -\off 
\sum_{\pff k\qff \in\qff C\vphantom{K^K}}\qff m_{\dff k}\qff.
\]

\vspace{-12pt} 
If\oss $i\qff \not\in\pff C$\nnsp,\qff\oss then\qss 
$m_{\dff i}\off =\off 0$\qss 
and\dss hence\vspace{3pt} 
\[
\quad
x_{\dff i}\qff -\qff m_{\dff i}
\off\off =\off\off
x_{\dff i}
\off\off \leq\off\off 
\sum_{k\qff \not\in\qff C\vphantom{K^K}}\qff x_{\dff k}
\off\off =\off\off
1\off -\off 
\sum_{k\qff \in\qff C\vphantom{K^K}}\qff x_{\dff k}
\off\off \leq\off\off 
1\off -\off 
\sum_{k\qff \in\qff C\vphantom{K^K}}\qff m_{\dff k}\qff.
\]

\vspace{-12pt} 
Together\dss with\qss (\ref{c-sum})\qss
these inequalities\sss imply\dss that\qss
$x_{\dff i}\qff -\qff m_{\dff i}\off <\off 
\num{C}\fff\left/\dff l\right.$\qss
for every\qss $i\qff \in\pff I$\nnsp.\oss  \eproof

\myuppar{Continuous self-maps of\dss $\Delta^n$\dnsp.}
Now\dss we\sss turn\dss to\sss the final\dss part\sss of\trs the proof\dss
of\qss Brouwer's\qss theorem.\oss
Let\dss
$f\dff \colon\dff 
\Delta^n\qff \ttoo\qff \Delta^n$\qss
be a continuous map.\oss
Recall\dss that\qss
$T\off =\off T_{\fff l}$\qss
depends on\dss $l$\nnsp.\qff\oss
If\vspace{3pt}
\[
\quad
x\off =\off
(\qff x_{\trf 0}\fff,\pff  x_{\dff 1}\fff,\pff \ldots\fff,\pff  x_{\dff n} \qff)
\qff \in\pff T_{\fff l}
\hspace*{1em}\mbox{and}\hspace*{1em}
y\off =\off
(\qff y_{\trf 0}\fff,\pff  y_{\dff 1}\fff,\pff \ldots\fff,\pff  y_{\dff n} \qff)
\off =\off
f\dff(\dff x\trf)
\qff,
\]

\vspace{-9pt}
then\qss
$x_{\trf 0}\qff +\qff  x_{\dff 1}\qff +\qff \ldots\qff +\qff  x_{\dff n}
\off =\off
y_{\trf 0}\qff +\qff  y_{\dff 1}\qff +\qff \ldots\qff +\qff  y_{\dff n}
\off =\off
1$\qss
and\dss hence\qss
$y_{\dff i}\off \geq\off x_{\dff i}$\qss for some\qss $i\qff \in\pff I$\nnsp.\oss
Let\sss $c\dff(\dff x\trf)$\dss be equal\dss to any\sss such $i$\nnsp.\oss
This rule defines\sss  
a coloring
of\trs $T_{\fff l}$\nsp,\oss i.e.\qss a map\qss 
$c\dff \colon\dff 
T_{\fff l}\qff \ttoo\qff I$\nnsp.\oss

By\qss Theorem\qss \ref{main-lemma}\qss for every\sss $l$\dss
there\sss exists\sss a\sss properly\sss colored\sss cell,\oss
i.e.\qss a\sss pair\dss
$(\trf \sigma_{\fff l}\dff,\qff C_{\dff l} \trf)$\dss
such\dss that\dss
$C_{\dff l}$\dss is\dss non-empty\halfff,\pss
$\sigma_{\fff l}\off \subset\off T_{\fff l}$\qss is\dss dominant\dss with\sss respect\dss to\dss
$C_{\dff l}$\nsp,\oss
and\qss
$C_{\dff l}\off =\off c\dff(\dff \sigma_{\fff l}\trf)$\nnsp.\oss
By\qss Theorem\qss \ref{dominant-size}\qss
the diameter of\trs the sets\dss $\sigma_{\fff l}$\dss
tends\sss to $0$ when\qss
$l\qff \toto\qff \infty$\nnsp.\oss
Therefore,\oss
after\dss passing to a subsequence,\oss still denoted\dss by\dss $\sigma_{\fff l}$\nnsp,\oss 
we can assume\sss that\sss all\sss elements of\dss
$\sigma_{\dff l}$ converge to the same point\qss 
$z
\off =\off
(\qff z_{\trf 0}\fff,\pff  z_{\dff 1}\fff,\pff \ldots\fff,\pff z_{\dff n} \qff)
\qff \in\pff \Delta^n$\qss when\qss 
$l\qff \toto\qff \infty$\nnsp.\oss
Let\vspace{3pt}
\[
\quad
w\off\dff =\off
(\qff w_{\trf 0}\fff,\pff  w_{\dff 1}\fff,\pff \ldots\fff,\pff  w_{\dff n} \qff)
\off =\off\dff
f\dff(\dff z\trf)
\qff.
\]

\vspace{-9pt}
Since there\sss are\sss only\dss  finitely\dss many\sss subsets of\dss $I$\nnsp,\oss
after\dss passing\dss to a\sss further subsequence we can assume\sss that\qss
$C_{\dff l}\off =\off C$\qss for some non-empty\sss subset\qss
$C\off \subset\off I$\qss
independent\sss of\trs $l$\nnsp.\oss
Then\qss
$C
\off =\off 
C_{\dff l}
\off =\off 
c\dff(\dff \sigma_{\fff l}\trf)$\qss
for every\sss $l$\dss
and\dss hence for every\sss $l$\dss and every\qss $i\qff \in\pff C$\qss
there\dss is\dss a\sss point\qss
$z\dff(\dff i\fff,\qff l\qff)\qff \in\pff \sigma_{\dff l}$\qss 
such\dss that\qss
$c\dff(\trf z\dff(\dff i\fff,\qff l\qff) \trf)\off =\off i$\nnsp.\oss
By\dss the choice of\trs the colorings $c$\sss
passing\dss to\sss the limit\qss
$l\qff \toto\qff \infty$\qss
shows\sss that\qss
$w_{\dff i}\off \geq\off z_{\dff i}$\qss
for every\qss $i\qff \in\pff C$\nnsp.\oss

At\dss the same\sss time\qss Theorem\qss \ref{coordinate-orders}\qss implies\sss that\dss
$x_{\dff i}\qff <\qff (\dff n\qff +\qff 1\dff)\fff/\fff l$\oss
for\dss every\qss $x\qff \in\qff \sigma$\qss 
and\qss $i\qff \in\pff I\qff \smallsetminus\qff C$\nnsp.\oss
By passing to the limit\qss 
$n\qff \toto\qff \infty$\qss 
we conclude\sss that\qss
$z_{\dff i}\qff =\qff 0$\qss 
for every\qss $i\qff \in\pff I\qff \smallsetminus\qff C$\nnsp.\oss
Therefore\vspace{1.5pt}
\[
\quad
\sum\nolimits_{\pff i\qff \in\qff C\vphantom{K^K}}\qff z_{\dff i}
\off\off =\off\off
1\dff.
\]

\vspace{-9pt}
Since\qss
$w_{\dff i}\off \geq\off z_{\dff i}$\qss
for every\qss $i\qff \in\qff C$\nnsp,\oss
this equality\dss implies\sss that\oss 
$\dis
\sum\nolimits_{\pff i\qff \in\qff C\vphantom{K^K}}\qff w_{\dff i}
\off\off \geq\off\off
1$\nnsp.\oss
But\vspace{0.25pt}
\[
\quad
\sum\nolimits_{\pff i\qff \in\qff I\vphantom{K^K}}\qff w_{\dff i}
\off\off =\off\off
1 
\qff.
\]

\vspace{-10.5pt}
It\dss follows\dss that\qss 
$w_{\dff i}\off =\off 0\off =\off z_{\dff i}$\qss
for every\qss $i\qff \in\pff I\qff \smallsetminus\qff C$\qss
and\oss 
$\dis
\sum\nolimits_{\pff i\qff \in\qff C\vphantom{K^K}}\qff w_{\dff i}
\off\off =\off\off
1
\off\off =\off\off
\sum\nolimits_{\pff i\qff \in\qff C\vphantom{K^K}}\qff z_{\dff i}$\nsp.\oss

Since\qss
$w_{\dff i}\off \geq\off z_{\dff i}$\qss
for every\qss $i\qff \in\qff C$\nnsp,\oss
this equality\dss implies\sss that\qss
$w_{\dff i}\off =\off z_{\dff i}$\qss
for every\qss $i\qff \in\pff C$\qss also.\oss
Therefore\qss
$w_{\dff i}\off =\off z_{\dff i}$\qss
for every\qss $i\qff \in\pff I$\qss
and\dss hence\qss
$f\dff(\dff z\trf)\off =\off w\off =\off z$\nnsp,\oss
i.e.\dss $z$\dss is\dss a\sss fixed\dss point\sss of\trs the map\dss $f$\nnsp.\oss
This completes\sss the proof\dss of\qss Brouwer's\qss fixed\dss point\dss theorem.\oss  \eproof

\mysection{Scarf's\qss proof\qss and\qss its\qss versions}{comparing}

\myuppar{Preferences and\dss utility\dss functions.}
In\dss the context\sss of\trs the mathematical\sss economics and\dss the game\sss theory\dss
the elements of\trs the set\sss $T$ can\dss be interpreted as goods,\oss
and\dss the elements of\trs the set $I$ as\sss traders on a market\halfff.\oss
From\dss this point\sss of\dss view,\oss the order\dss $<_{\dff i}$\dss on\dss the set\sss
of\dss goods reflects\sss the preferences of\trs the\sss trader\dss $i$\nnsp.\oss
It\dss is\dss quite natural\dss to\sss think\dss that\dss the preferences of\dss 
each\dss trader\sss $i$
arise from a\dss \emph{utility\dss function}\qss
$u_{\dff i}\dff \colon\dff
T\qff \ttoo\qff 
\rrr_{\qff +}
\off =\off 
(\dff 0\fff,\pff \infty\dff)$\dss
such\dss that\dss the inequalities\qss $x\qff <_{\dff i}\qff y$\qss
and\qss $u_{\dff i}\dff(\dff x\trf)\qff <\qff u_{\dff i}\dff(\dff y\trf)$\qss
are equivalent\halfff.\oss
In any\sss case,\oss such\sss 
utility\dss functions always exist\halfff.\oss
By\sss numbering\dss the\sss traders 
one can arrange utility\dss functions into a map\qss
$u\dff \colon\dff
T\qff \ttoo \rrr^{\fff n\dff +\dff 1}_{\qff +}$\qss
taking\dss $x\qff \in\pff T$\dss to\sss the point\trs  
$u\dff(\dff x\trf)\qff \in\qff \rrr^{\fff n\dff +\dff 1}_{\qff +}$\dss
having\dss the utilities $u_{\dff i}\dff(\dff x\trf)$ 
as\sss its coordinates.\oss
By\dss perturbing\sss a\sss little\sss the utility\dss functions,\oss
if\trs necessary\halfff,\oss one can assume\sss that\sss $u$
is\dss injective,\oss and\dss then\dss identify\dss the set\sss $T$
with\dss its\dss image $u\dff(\trf T\trf)$\nnsp.\oss
Therefore,\oss without\sss any\dss loss of\dss generality\sss one can
assume\sss that\sss $T$ is\dss a\sss finite subset\sss of\dss 
$\rrr^{\fff n\dff +\dff 1}_{\qff +}$\dss and\dss the orders\dss $<_{\dff i}$\dss
are determined\dss by\dss the values of\dss coordinates.\oss 
For applications\sss to\qss Brouwer's\qss fixed\dss point\trs theorem\dss
the case of\qss $T\off \subset\off \Delta^n$\qss is\dss sufficient\halfff.\oss
This\dss is\dss the framework of\qss Scarf's\qss papers\qss
\cite{sc1}\qss --\qss \cite{sc3}.\oss

\myuppar{Slack\dss vectors.}
Scarf\qss enlarges\sss $T$\sss
by\dss adding\dss to\sss it\dss vectors representing\dss the
$(\fff n\dff -\dff 1\fff)$\dnsp-faces of\dss $\Delta^n$\dnsp.\oss
The face defined\dss by\dss the equation\qss
$x_{\dff i}\off =\off 0$\qss is\dss represented\dss by\vspace{3.75pt}
\[
\quad
s\dff(\dff i\trf)
\off =\off
(\qff 
M_{\dff i}\fff,\pff 
\ldots\fff,\pff
M_{\dff i}\fff,\pff
0\fff,\pff
M_{\dff i}\fff,\pff
\ldots\fff,\pff
M_{\dff i}
\qff)
\qff \in\pff \rrr^{\fff n\dff +\dff 1}
\qff,
\]

\vspace{-8.25pt}
the vector\dss with\dss the $i${\dnsp}th\dss coordinate equal\dss to $0$
and other coordinates equal\dss to some real\dss number\dss $M_{\dff i}$\dss
bigger\dss than\dss the $i${\dnsp}th\dss coordinate of\dss every\dss vector\qss
$x\qff \in\pff T$\dnsp.\oss
If\qss $T\off \subset\off \Delta^n$\nnsp,\oss
it\dss is\dss sufficient\dss to assume\sss that\qss 
$M_{\dff i}\qff >\qff 1$\nnsp.\oss
To stress an analogy\dss with\dss 
linear\dss programming\halfff,\oss Scarf\qss
calls\sss these vectors\qss \emph{slack\sss vectors}.\oss
Let\dss $\mathbb{I}$\dss be\sss the set\sss of\trs the slack\sss vectors.\oss

It\dss is\dss convenient\dss to identify\dss the $i${\nnsp}th\dss trader\dss
with\dss the $i${\nnsp}th\dss slack\sss vector\sss $s\trf(\dff i\trf)$\sss
for every\dss $i\qff \in\pff I$\dss and\dss therefore\dss identify\sss $I$\sss
with\dss $\mathbb{I}$\nnsp.\oss
The en\-larged\sss set\dss $T$\dss is\dss  
the union\qss
$T\qff \cup\pff \mathbb{I}$\nnsp.\oss
If\trs  
$M_{\dff i}$\sss are pairwise different\halfff,\oss
then\dss the values of\dss coordinates define\sss linear\sss orders on\qss
$T\qff \cup\pff \mathbb{I}$\nnsp.\qff\oss
A\sss subset\pss 
$X\off \subset\off T\qff \cup\pff \mathbb{I}$\pss is\dss called\qss
\emph{primitive}\oss if\pss 
$\num{X}\off =\off n\qff +\qff 1$\pss and\dss
$X$\dss is\dss a\sss dominant\sss subset\sss of\qss
$T\qff \cup\pff \mathbb{I}$\dss 
with\dss respect\dss to\sss the set\dss $I\off =\off \mathbb{I}$\dss of\dss all\sss orders.

By\dss the choice of\dss slack\sss vectors 
a\sss subset\qss $\sigma\off \subset\off T$\qss
is\dss dominant\dss with\dss respect\dss to\qss
$C\off \subset\off \mathbb{I}$\qss
if\trs and\dss only\trs if\qss
$\sigma\qff \cup\qff (\qff \mathbb{I}\qff \smallsetminus\qff C\qff)$\qss
is\dss a\sss primitive set\halfff.\oss
Conversely\halfff,\pss
$X\off \subset\off T\qff \cup\pff \mathbb{I}$\qss
is\dss a\sss primitive subset\trs 
if\trs and\dss only\trs if\trs the intersection\qss
$X\qff \cap\qff T$\qss is\dss dominant\dss with\dss respect\dss to\qss
$\mathbb{I}\qff \smallsetminus\qff X$\nnsp.\oss
Therefore\sss the notions of\trs primitive and dominant\sss sets
are equivalent\sss and\sss one can use either\sss of\trs them.\oss

\myuppar{Scarf's\qss main\dss lemma.}
\emph{Suppose\sss that\qss
$X\off \subset\off T\qff \cup\pff \mathbb{I}$\pss
is\dss a\sss primitive set\halfff,\oss
and\qss $x\qff \in\pff X$\nnsp.\oss
Then either\qss
$X\qff -\qff x\off \subset\off \mathbb{I}$\nnsp,\oss
or\dss there\dss is\dss a\sss unique\qss
$y\qff \in\pff (\qff T\qff \cup\pff \mathbb{I} \qff)
\off \smallsetminus\off X$\qss such\dss that\qss
$X\qff -\qff x\qff +\qff y$\qss
is\dss a\sss primitive set\halfff.\oss
For every\qss $i\qff \in\pff \mathbb{I}$\dss
there\dss is\dss exactly\sss one\sss primitive set\sss
containing\dss $\mathbb{I}\qff -\qff i$\nnsp.\oss}

\vspace{7pt}
See\qss \cite{sc2},\oss Lemma\qss 1.\oss
The analogy\dss between\dss primitive sets and $n$\dnsp-simplices of\dss
a\sss triangulation of\dss $\Delta^n$\sss suggests\sss 
to state\sss this\sss lemma differently\halfff.\oss
Let\dss us\dss call\sss a\sss subset\sss
$Y\off \subset\off T\qff \cup\pff \mathbb{I}$\qss
\emph{almost\dss primitive}\oss
if\qss
$\num{Y}\off =\off n$\dss
and\dss there\dss is\dss a\sss primitive set\sss containing\sss $Y$\nnsp.\oss
The set\dss
$Y\off =\off X\qff -\qff x$\dss
from\qss Scarf's\qss main\dss lemma\dss is\dss almost\dss primitive,\oss
and\dss hence we can\dss restate\sss this\sss lemma as follows.\oss\vspace{1pt}

\myuppar{Another\dss form\sss of\pss Scarf's\qss main\dss lemma.}
\emph{Suppose\sss that\qss 
$Y\qff \subset\pff T\qff \cup\pff \mathbb{I}$\dss
is\dss an\sss almost\dss primitive set\halfff.\oss
Then either\oss $Y\off \subset\off \mathbb{I}$\nnsp,\oss
or\dss $Y$\dss is\dss a subset\dss of\qss exactly\dss
two primitive sets.\oss
For every\qss $i\qff \in\pff \mathbb{I}$\dss
the set\qss $\mathbb{I}\qff -\qff i$\qss
is\dss almost\dss primitive and\dss
is\dss a subset\dss of\qss only\sss one primitive set\halfff.\oss}

\vspace{7pt}
By\dss the definitions,\oss 
if\qss 
$X\off \subset\off T\qff \cup\pff \mathbb{I}$\qss 
is\dss a\sss primitive set\halfff,\oss
then\dss 
$(\qff X\qff \cap\qff T\fff,\pff \mathbb{I}\qff \smallsetminus\qff X\qff)$\dss
is\dss a\sss room,\oss
and\dss if\qss $X$\dss is\dss an\sss almost\dss primitive set\halfff,\oss
then\sss this\dss pair\dss is\dss a\sss door\halfff.\oss
Moreover\halfff,\oss if\trs $Y$\sss is\dss an almost\dss primitive set\halfff,\pss
$X$\sss is\dss a\sss primitive set\halfff,\oss 
and\qss $Y\off \subset\off X$\nnsp,\oss
then\dss  
$(\qff Y\qff \cap\qff T\fff,\pff \mathbb{I}\qff \smallsetminus\qff Y\qff)$\dss
is\dss a\sss door of\trs the room\dss
$(\qff X\qff \cap\qff T\fff,\pff \mathbb{I}\qff \smallsetminus\qff X\qff)$\nnsp.\oss
Conversely\halfff,\oss if\qss
$(\dff \sigma\fff,\qff C \dff)$ is\dss a\sss room,\oss
then\qss $\sigma\qff \cup\qff (\qff \mathbb{I}\qff \smallsetminus\qff C\qff)$\qss
is\dss a\sss primitive set\halfff.\oss
The corresponding statement\dss for\dss the doors\dss is\dss only\sss
a\sss little\sss less obvious.\oss
By\qss Lemmas\qss \ref{outside-doors}\qss and\qss \ref{doors}\qss
every\sss door\dss is\dss a\sss door\sss of\dss a\sss room\qss
(one can also see\sss this directly\fff).\oss
It\dss follows\dss that\dss if\qss
$(\dff \tau\fff,\qff D \dff)$ is\dss a\sss door\halfff,\oss
then\qss $\tau\qff \cup\qff (\qff \mathbb{I}\qff \smallsetminus\qff D\trf)$\qss
is\dss contained\dss in\sss a\sss primitive set\sss and\dss hence\dss is\dss
an\sss almost\sss primitive set\halfff.\oss
So,\oss
the notions of\trs primitive and\sss almost\dss primitive sets\sss
are equivalent\dss to\sss 
the notions of\trs rooms and\sss doors,\oss and\qss
Scarf's\qss main\dss lemma\dss
is\dss equivalent\dss to\qss 
Lemmas\qss \ref{outside-doors}\qss and\qss \ref{doors}\qss 
together\halfff.\oss\vspace{1.75pt}

\myuppar{Primitive\sss sets and\dss rooms.}
An obvious advantage of\dss working\dss with\sss slack\sss vectors\sss
and\dss primitive sets\dss is\dss
the fact\dss that\sss 
an almost\dss primitive set\dss is\dss always 
obtained\dss by\dss removing a point\dss
from a primitive set\halfff.\oss
In contrast\halfff,\oss a door\dss can\dss be obtained\dss from a room
$(\dff \sigma\fff,\qff C \dff)$ 
either\dss by\dss removing a point\dss from $\sigma$ 
or\trs by\sss adding\dss a point\dss to $C$\nnsp.\oss
Still\halfff,\oss traders and\dss goods are\sss inherently\sss different\halfff,\pss
and\dss it\dss is\dss only\dss natural\dss to keep\sss track of\trs this distinction,\oss
as\dss it\dss is\dss done in\dss the notion of\dss a\sss room.\oss

Scarf's\qss proof\dss of\trs his combinatorial\dss theorem\dss
uncovered a structure\sss hidden\dss in a\sss family\sss 
of\dss orders on a\sss given set,\oss
namely\halfff,\oss the collection of\dss primitive sets.\oss
H.\qff~Kuhn\qss \cite{ku}\qss interpreted\trs this structure
as\sss the structure of\dss a\sss simplicial\sss complex\halfff.\oss
The $n$\dnsp-simplices of\trs this complex\sss 
are primitive sets 
and\qss Scarf's\qss main\dss lemma\sss means\sss
that\dss this complex\dss is\dss a\sss pseudo-manifold.\oss
Sperner's\qss lemma\dss and\dss its proofs naturally\sss extend\dss
to\sss pseudomanifolds\dss 
and\qss Scarf's\qss combinatorial\dss theorem\dss
can\dss be proved\dss by\sss applying\dss such an extension\dss
to\sss this pseudo-manifold.\oss

By\dss keeping\dss the distinction\dss between\dss goods and\dss traders one
can see an even\dss more rich combinatorial\sss structure\sss behind\dss the
family\sss of\dss orders.\oss
Namely\halfff,\oss a\sss room $(\dff \sigma\fff,\qff C \dff)$
can\dss be\sss thought\sss as\qss ``belonging''\qss to\sss the face $C$
of\trs the abstract\sss simplex\dss having\dss $\mathbb{I}$\dss as\sss the
set\sss of\dss vertices.\oss
This\sss leads\sss to a\sss refined\sss version of\pss Kuhn's\qss interpretation
and\dss to a\sss proof\dss of\pss Scarf's\qss combinatorial\trs theorem\dss
using\qss Scarf's\qss main\dss lemma as\sss the starting\dss point\sss
and\dss then\dss proceeding\dss in\dss the same way as\sss in\sss
combinatorial\sss or\dss homological\dss proofs\sss of\pss
Sperner's\qss lemma.\oss
See\qss \cite{i2},\oss Section\qss 1.\oss

\myuppar{Petri--Voorneveld\qss version.}
In contrast\dss with\qss Scarf\qss \cite{sc2}\qss
and\dss with\dss Section\qss \ref{scarf},\oss
H.\dss Petri\dss and\dss M.\dss Voorneveld\qss \cite{pv}\qss
start\dss with\dss fixing\sss a coloring\qss
$c\dff \colon\dff
T\qff \ttoo\qff I$\nnsp.\oss
In\dss the language of\trs the present\dss paper\halfff,\oss
they\dss work only\dss with\dss pairs $(\dff \sigma\fff,\qff C \dff)$
such\dss that\qss
$\num{C\qff \smallsetminus\qff c\dff(\trf \sigma \trf)}\off \leq\off 1$\nnsp.\oss
Petri\dss and\dss Voorneveld\dss call\trs them\qss
\emph{candidates}.\oss
One can easily\dss check\dss that\sss every\sss such\dss pair\dss is\dss
either a room or a door\halfff,\oss
but\dss other\dss rooms and\sss doors are missing\halfff.\oss
Of\dss course,\oss such\dss pairs\sss are\sss sufficient\dss for\dss
proving\qss Scarf\qss combinatorial\dss theorem\qss
(called\dss in\pss \cite{pv}\pss \emph{No-bullying\dss lemma}\halfff\fff),\oss
but\dss the structure uncovered\dss by\qss Scarf\qss remains\sss hidden.\oss
Perhaps,\oss this reflects\sss the intention\dss to eliminate geometry\halfff.\oss

One has\sss to admit\dss that\dss the previous paragraph\dss is\dss not\dss
quite faithful\dss to\qss \cite{pv}.\oss
After\dss fixing\sss a coloring\halfff,\oss
Petri\dss and\dss Voorneveld\dss replace each color\qss
(i.e.\qss each of\trs the orders)\qss
by\sss as many\sss of\trs its\sss copies as\sss there are elements of\trs $T$\sss 
with\dss this color\halfff.\oss
This\sss trick\sss allows\sss to\sss turn\dss the coloring $c$ into
a\sss bijection.\oss
Petri\dss and\dss Voorneveld\dss go one step\sss further and\dss use\sss
this\sss bijective coloring\dss to identify\trs $T$\sss and\qss $I$\nnsp.\oss
This erases\sss the difference between\dss goods and\dss traders completely\qss
(actually\halfff,\oss Petri\dss and\dss Voorneveld\dss speak about\dss toys and\dss kids).\oss
For\dss the present\sss author\dss this\sss identification\dss was\sss the main\sss
stumbling\dss block\dss in\dss reading\qss \cite{pv}.\oss

\myuppar{The\sss paths\dss to\sss primitive sets and colorful\dss rooms.}
Let\dss
$c\dff \colon\dff
T\qff \ttoo\qff \mathbb{I}$\dss
be\sss a coloring\halfff,\oss
and\dss let\dss us extend $c$\sss to a coloring\dss
$T\qff \cup\pff \mathbb{I}\qff \ttoo\qff \mathbb{I}$\dss
equal\dss to\sss the identity\sss on\dss $\mathbb{I}$\nnsp.\oss
There\dss is\dss no danger\dss in\sss denoting\dss this extension also by\sss $c$\nnsp.\oss
In\dss this\sss language\qss Scarf's\qss combinatorial\dss theorem
asserts\sss that\dss there\dss is\dss a\sss primitive set\dss $X$\dss
such\dss that\dss $c\trf(\trf X\trf)\off =\off \mathbb{I}$\nnsp.\oss
In order\dss to prove\sss this,\oss
Scarf's\qss uses as\sss the main\dss tool\dss the operation of\dss
replacing\dss $X$\dss by\qss
$X\qff -\qff x\qff +\qff y$\qss as\sss in\dss the main\dss lemma.\oss
This\sss replaces\sss in\dss $X$\dss the element\sss $x$\dss by\sss
a new element\sss $y$\sss
and\dss is\dss called\sss a\qss
\emph{replacement\dss step}.\oss 
The proof\dss starts with\dss fixing
a color\dss $i\qff \in\pff \mathbb{I}$\dss
and considering\dss the set\dss $\mathbb{I}\qff -\qff i$\nnsp.\oss 
Clearly\halfff,\oss the set\dss $\mathbb{I}\qff -\qff i$\dss is\dss almost\dss primitive and\dss
by\dss the main\dss lemma\dss is\dss contained\dss
in a unique primitive set\dss $X$\nnsp.\oss
If\qss $c\trf(\trf X\trf)\off =\off I$\nnsp,\oss then we are done.\oss
Otherwise\qss
$c\trf(\trf X\trf)\off =\off \mathbb{I}\qff -\qff i$\dss
and\sss one can apply\dss to\dss
$X$\dss a sequence of\trs
replacement\sss steps.\oss
The elements\sss to replace are chosen\dss in such a way\dss that\sss
after each\sss step either\qss $c\trf(\trf X\trf)\off =\off I$\qss
or\qss $c\trf(\trf X\trf)\off =\off \mathbb{I}\qff -\qff i$\nnsp.\oss
Scarf\qss proves\sss that\dss such sequence of\trs replacement\sss steps 
eventually\dss
arrives\sss at\sss a\sss primitive set\dss $X$\dss such\dss that\qss
$c\trf(\trf X\trf)\off =\off I$\nnsp,\oss
when\dss it\sss stops.\oss

One can split\dss each\dss replacement\sss step\sss into\sss two
simpler steps.\oss
Namely\halfff,\oss
instead\sss of\dss passing\dss from a\sss primitive set\dss $X$\dss
to\dss $X\qff -\qff x\qff +\qff y$\dss
as in\dss the\dss main\dss lemma,\oss
one can\dss first\dss pass from\dss $X$\dss to an almost\dss primitive set\qss
$X\qff -\qff x$\dss and\dss then\sss pass from\dss $X\qff -\qff x$\qss
to\qss $X\qff -\qff x\qff +\qff y$\nnsp.\oss
Let\dss us\dss split\dss in\dss this way\sss every\dss replacement\sss step in\qss
Scarf's\qss sequence and add\qss $\mathbb{I}\qff -\qff i$\qss at\dss the beginning\halfff.\oss
The resulting\sss sequence 
alternates between almost\dss primitive and\dss primitive sets and\qss
$c\trf(\trf X\trf)\off =\off\qff \mathbb{I}\qff -\qff i$\qss for every\dss its\dss term\dss $X$\dss
except\sss of\trs the\sss last\sss one.\oss
The transformation\qss
$Z\off \longmapsto\off
(\qff Z\qff \cap\qff T\fff,\pff \mathbb{I}\qff \smallsetminus\qff Z\qff)$\qss
turns\sss this sequence\qss
into a sequence rooms and\sss doors
starting\dss with\dss 
$(\trf \varnothing\fff,\pff  \{\trf i \qff\} \trf)$\dss
and ending\dss with a colorful\dss room.\oss
The condition\qss
$c\trf(\qff X\trf)
\off =\off\qff
\mathbb{I}\qff -\qff i$\qss
implies\sss that\dss doors and\dss rooms\dss
in\dss this sequence\dss are\dss nearly\sss colorful\sss
of\trs type $i$\nnsp,\oss except\dss the last\sss one.\oss
Moreover\halfff,\oss
this sequence\dss 
is\dss a\sss path\dss in\dss the graph\dss $G_{\dff i}$\dss
connecting\dss $(\trf \varnothing\fff,\pff  \{\trf i \qff\} \trf)$\dss
with a colorful\dss room.\oss
But\qss Theorem\qss \ref{graph}\pss implies\sss that\dss there\dss is\dss only\sss
one\sss such\dss path.\oss 
Hence\sss the\sss transformation\qss
$Z\off \longmapsto\off
(\qff Z\qff \cap\qff T\fff,\pff \mathbb{I}\qff \smallsetminus\qff Z\qff)$\qss
turns\qss Scarf's\qss sequence\sss
into\dss the\sss path\dss to a colorful\dss room\dss
from\dss Section\qss \ref{scarf}.\oss
In other\dss words,\oss
the arguments of\qss Section\qss \ref{scarf}\qss
lead\dss to essentially\dss the same\sss the path\dss to a colorful\dss room\sss
as\qss Scarf's\qss ones.\oss
Of\dss course,\oss the same\dss is\dss true with\dss respect\dss to\sss
the arguments of\pss
Petri\dss and\dss Voorneveld\pss \cite{pv}.\oss

\vspace*{12pt}
\begin{flushright}

August\qss 17,\oss 2019.\off\oss
Preface:\oss July\qss 21,\oss 2022
 
https\halfff:/\!/\hspace*{-0.06em}nikolaivivanov.com

E-mail\halfff:\oss nikolai.v.ivanov{\fff}@{\dff}icloud.com,\oss ivanov{\fff}@{\dff}msu.edu

Department\sss of\qss Mathematics,\oss Michigan\sss State\sss University

\end{flushright}

\end{document}